\documentclass[10pt]{article}
\date{24.11.2014} 
\usepackage{hyperref}
\usepackage{amsfonts,amssymb,amsmath,amscd,euscript,array,mathrsfs}
\input{liemacs10.sty} 
\renewcommand{\oline}{\overline} 
 
\newcommand{\cJ}{\mathcal J} 

\newcommand{\bD}{\mathbb D}

\newcommand{\vphi}{\varphi} 
\renewcommand{\phi}{\varphi} % to be removed in the end 
\renewcommand{\subeq}{\subseteq} 
\renewcommand{\la}{\langle}
\renewcommand{\ra}{\rangle}

\renewcommand{\mlabel}{\label}

\newcommand\CP{\mathop{\rm CP}\nolimits}
\newcommand\ball{\mathop{\tt ball}\nolimits}
\renewcommand\P{\mathop{\rm P}\nolimits}

\begin{document} 

\title{Nonlinear completely positive maps and dilation theory \\ for  real 
involutive algebras} 
\author{Daniel Belti\c t\u a and Karl-Hermann Neeb} 

\maketitle 
%%%%%%%%%%%%%%%%%%%%%%%%%%

\begin{abstract} A real seminormed involutive algebra is a 
real associative algebra $\cA$ endowed with an involutive antiautomorphism 
$*$ and a submultiplicative seminorm $p$ with $p(a^*) =p(a)$ for $a\in \cA$. 
Then $\ball(\cA,p) := \{ a \in \cA \: p(a) < 1\}$ is an involutive subsemigroup. 
For the case where $\cA$ is unital, our main result asserts that a function 
$\phi \: \ball(\cA,p) \to B(V)$, $V$ a Hilbert space, is completely 
positive (defined suitably) if and only if it is positive definite and 
analytic for any locally convex topology for which $\ball(\cA,p)$ is open. 
If $\eta_\cA \: \cA \to C^*(\cA,p)$ is the enveloping $C^*$-algebra of $(\cA,p)$ 
and $e^{C^*(\cA,p)}$ is the $c_0$-direct sum of the symmetric tensor powers 
$S^n(C^*(\cA,p))$, then the above two properties are equivalent to the 
existence of a factorization $\phi = \Phi \circ \Gamma$, 
where $\Phi \: e^{C^*(\cA,p)} \to B(V)$ is linear completely positive and 
$\Gamma(a) = \sum_{n = 0}^\infty \eta_\cA(a)^{\otimes n}$. We also obtain a suitable generalization 
to non-unital algebras. 

An important consequence of this result is a description of the  
unitary representations of $\U(\cA)$ with 
bounded analytic extensions to $\ball(\cA,p)$ in 
terms of representations of the $C^*$-algebra $e^{C^*(\cA,p)}$. 
\end{abstract} 

\tableofcontents

\section{Introduction} \mlabel{sec:intro}
%\todo{to be written} 

A real {\it seminormed involutive algebra (or $*$-algebra)} 
is a pair $(\cA,p)$, consisting 
of a real associative algebra $\cA$ endowed with an involutive antiautomorphism 
$*$ and a submultiplicative seminorm $p$ satisfying $p(a^*) =p(a)$ for $a\in \cA$. 
Our project eventually aims at a systematic understanding of unitary 
representations of unitary groups $\U(\cA)$ of real seminormed involutive algebras 
$(\cA,p)$. 
If $\cA$ is unital, its unitary group is 
\[ \U(\cA) := \{ a \in \cA : a^* a = aa^* = \1\}, \] 
and if $\cA$ is non-unital, then $\U(\cA):=\U(\cA^1)\cap(\1+\cA)$, 
where $\cA^1 = \cA \oplus \R \1$ is the unitization of $\cA$.   
Typical examples we have in mind are  $C^*$-algebras (considered as real algebras) 
and algebras of the form 
$\cA=C^\infty(X,\cB)$ for a Banach $*$-algebra $\cB$ or  
$\cA = C^\infty(X,M_n(\K))$, where $X$ is a smooth manifold and 
$\K \in \{\R,\C,\H\}$. 
In the latter case $\U(\cA) \cong C^\infty(X,\U_n(\K))$ are groups of smooth 
maps with values in a compact Lie group, hence particular examples of gauge groups. 

We wish to initiate a line of investigation of unitary representations 
of $\U(\cA)$ which, in a suitable sense, 
are boundary values of representations of the ball semigroups 
$\ball(\cA) := \{ a \in \cA \: p(a) < 1\}$. This leads us naturally to 
completely positive maps on these semigroups. 
The present paper provides the foundation of that study, by developing the 
corresponding dilation theory and reducing it completely to the $C^*$-context. 
In a sequel we will apply these dilation methods to 
the representation theory of semigroups $\ball(\cA,p)$ for 
various types of concrete algebras~$\cA$. 

We now introduce some concepts needed to state our main theorem. 
If $S$ is an involutive semigroup and $V$ a Hilbert space, then a 
positive definite function $\phi \: S \to B(V)$ is called {\it dilatable} 
if there exist a representation $(\pi,\cH)$ of $S$ and a continuous linear map 
$\iota \: V \to \cH$ such that $\phi(s) = \iota^* \pi(s) \iota$ for $s \in S$. 
We write $\eta_\cA \: \cA \to C^*(\cA,p)$ 
for the enveloping $C^*$-algebra of the real seminormed $*$-algebra 
$(\cA,p)$, and for a $C^*$-algebra $\cB$, we write $e^{\cB}$ for the $c_0$-direct sum of the $C^*$-algebras 
$S^n(\cB) \subeq \cB^{\otimes n}$, where the tensor products are constructed from 
the maximal $C^*$-cross norm (see \cite{Arv87} and Definition~\ref{exp_alg} below). 

Here is our main theorem: 
\begin{thm}\label{thm:intro} 
Let $(\cA,p)$ be a real seminormed involutive algebra and $V$ a complex Hilbert space.  
For a bounded function $\phi \: \ball(\cA,p) \to B(V)$,  the following are equivalent: 
\begin{itemize}
\item[\rm(i)] $\phi$ is completely positive and dilatable. 
\item[\rm(ii)] $\phi$ is dilatable, positive definite and analytic 
with respect to some locally convex topology for which $p$ is continuous. 
\item[\rm(iii)] There exists a linear completely positive map 
$\Phi \: e^{C^*(\cA,p)}\to B(V)$ with $\Phi \circ \Gamma = \phi$, 
where $\Gamma(a) = \sum_{n=0}^\infty \eta_\cA(a)^{\otimes n}$. 
\end{itemize}
If $\cA$ is unital, then every bounded positive definite function on $\ball(\cA,p)$ 
is dilatable, so that this requirement can be omitted in {\rm(i)} and {\rm(ii)}. 
\end{thm}

Note that (ii) is particularly strong because analyticity is not required for 
the topology defined by $p$ but for some topology which may be much finer. 
An important consequence of this theorem is that it leads to a one-to-one 
correspondence between bounded analytic representations of 
$\ball(\cA,p)$ with the representations of the $C^*$-algebra 
$e^{C^*(\cA,p)}$ (Theorem~\ref{final0}), which thus plays the role of a host algebra 
in the sense of \cite{Gr05} for the bounded analytic representations 
of the involutive semigroup $\ball(\cA,p)$, resp., the corresponding unitary 
representations of $\U(\cA)$, obtained from its action on $\ball(\cA,p)$ by multipliers.

To see how such representations arise naturally, consider 
the unitary group $\U(\cA)$ of a $C^*$-algebra $\cA$. Then 
every irreducible representation $(\pi, \cH)$ of $\cA$ provides 
for every partition $\lambda = (\lambda_1, \ldots, \lambda_k)$ of $N \in \N$ 
an irreducible representation $(\pi_\lambda, \bS_\lambda(\cH))$ on a subspace 
of $\cH^{\otimes N}$ by a straightforward generalization 
of the classical Schur--Weyl theory 
to Hilbert spaces (see \cite{BN12}, and \cite{EI13, Nes13} for an extension 
to type II$_1$-factor representations).  
Any such representation extends to a (holomorphic) polynomial 
representation of the multiplicative semigroup $(\cA,\cdot)$, resp., 
the $C^*$-algebra $S^N(\cA)$ on $\bS_\lambda(\cH)$. 

In the special case where 
$\cA = K(\cH)$ is the $C^*$-algebra of compact operators on 
an infinite dimensional complex Hilbert space $\cH$, 
results of A.~Kirillov and G.~Olshanski \cite{Ki73, Ol78} 
provide a complete classification of all unitary representation 
of $\U(\cA)$. All continuous unitary 
representations of $\U(\cA)$ are direct sums of irreducible ones 
which are of the form  \break $\bS_\lambda(\cH) \otimes \bS_\mu(\cH)^*$ 
and there are natural generalizations to $\U(K(\cH))$, where
$\cH$ is a real or a quaternionic Hilbert space. Therefore the 
Schur--Weyl construction is exhaustive in these three cases. 
The key method to obtain these results is to show that continuous unitary 
representations of $\U(K(\cH))$ are generated by dilation 
from operator-valued positive definite 
functions $\phi$ of the form $\phi(g) = \rho(ege)$, where
$e$ is a hermitian  projection of finite rank and 
$\rho$ is a continuous representation of the semigroup 
$S = \overline{\ball}(e\cA e)$ of contractions in $e\cA e$, 
where $e\cA e \cong M(n,\K)$ for $\K \in \{\R,\C,\H\}$. Comparing 
with the results of the present paper, it follows that 
the representations $\rho$ of $S$ for which $\phi$ is positive 
definite are precisely the completely positive ones 
(see \cite{Ne13} for a recent presentation of the Olshanski--Kirillov theory from this 
point of view). We conclude that, in this case, 
$e^{C^*(K(\cH))} \cong e^{K(\cH_\C)}$ is a host algebra for all continuous 
unitary representations of $\U(K(\cH))$. 

For the algebra $\cA = C^\infty(X, M_n(\K))$, $X$ a compact smooth manifold, 
the identity component of $\U(\cA)$ is a product of an abelian group and the group $C^\infty(X, \SU_N(\K))_0$ whose norm continuous unitary representations have recently been classified (\cite{JN13}). 
Here the irreducible ones are finite tensor products of so-called 
evaluation representations, hence obtained by restricting  finite tensor products of 
irreducible algebra representations. 
In \cite{BN15}, we will see how all these particular results fit naturally 
into our present framework and how the representations of $S^n(\cB)$ 
for a $C^*$-algebra $\cB$ can be described in terms of representations of $\cB$.\\

For the proof of our main theorem, we need several results that are 
interesting in their own right: 
\begin{itemize}
\item Factorization of $*$-representations of degree $n$ of $(\cA,\cdot)$ 
bounded on $\ball(\cA,p)$ through a linear representations 
of the $C^*$-algebra $S^n(C^*(\cA,p))$. So the $C^*$-functor is compatible with 
this kind of non-linearity. 
\item If $\phi$ is completely positive and homogeneous of degree $\alpha \in \R$, then 
$\alpha \in \N_0$ and $\phi$ is polynomial. This is based on Arveson's method of 
iterated differences (\cite{Arv87}) and the case $\cA = \R$, for which 
the completely positive maps  
$\phi \: \oline{\ball}(\cA) = [-1,1] \to \R$ are determined in \cite{CR78}. 
\end{itemize}

Theorem~\ref{thm:intro} generalizes Arveson's work on non-linear states 
on balls of unital $C^*$-algebras (\cite{Arv87}) which we found very inspiring. 
Actually he uses a notion of complete positivity that is weaker than ours 
and does not imply positive definiteness although he uses it at some point in 
his arguments. We discuss this issue in Section~\ref{sec:7}. His result asserts 
that completely positive states $\phi \: \ball(\cA,p) \to \C$ factor through 
linear states of the $C^*$-algebra $e^{\cA} \otimes e^{\oline\cA}$. This 
mysterious occurrence of the $C^*$-algebra $\overline{\cA}$ 
is clarified by our present approach via real algebras in 
which $C^*(\cA) \cong \cA_\C \cong \cA \oplus \oline\cA$, so that 
$e^{C^*(\cA)} \cong e^{\cA} \otimes e^{\oline\cA}$ provides the bridge to our 
Theorem~\ref{thm:intro}(iii). 

Theorem~\ref{thm:intro} also extends work of Ando/Choi on  expanding nonlinear 
completely positive maps  $\phi \: \cA \to \cB$ between 
$C^*$-algebras as series $\phi = \sum_{mn,= 0}^\infty \phi_{m,n}$, 
where the maps $\phi_{m,n}$ are homogeneous polynomial and completely  positive. 
For the special case case where $\cA$ is a $C^*$-algebra, 
the equivalence of (i) and (iii) in Theorem~\ref{thm:intro} is due to 
Hiai/Nakamura (\cite{HN87}). If $\cA$ is a $C^*$-algebra, then the situation simplifies considerably 
because an application of Arveson's Extension Theorem (\cite[Thm.~1.2.3]{Arv69}) 
shows that every bounded completely positive map $\phi$ on $\ball(\cA,p)$ is dilatable. 

The structure of this paper is as follows. 
In Section~\ref{sec:2} we introduce the basic concepts, 
such as real seminormed involutive algebras, their ball semigroups 
and completely positive functions thereon. Section~\ref{sec:3} develops 
dilation theory on the abstract level of 
structured $*$-semigroups. These are involutive semigroups 
$S$ with a homomorphism $\gamma \: (0,1] \to S$, so that 
$S^\circ := \bigcup_{0 < r < 1} \gamma(r)S$ is a semigroup ideal. 
Here our main result is Theorem~\ref{thm:Arv2.2} on 
the existence of dilations for positive definite 
functions on $S^\circ$. This 
applies in particular to the subsemigroup $S^\circ = \ball(\cA,p)$ for the 
structured $*$-semigroup $S = \oline\ball(\cA,p)$, when $\cA$ is unital. 
In Section~\ref{sec:4} we  
provide the key tools to obtain a series expansion $\phi = \sum\limits_{n = 0}^\infty \phi_n$ 
of a bounded completely positive function $\phi$ on $\ball(\cA,p)$ into homogeneous 
components $\phi_n$. 

In Section~\ref{sec:5} we turn to positive definite functions 
$\phi \: \ball(\cA,p) \to B(V)$ which are analytic with respect to some 
locally convex topology on $\cA$ for which $\ball(\cA,p)$ is open. 
We derive a series expansion of dilatable bounded analytic positive definite functions, 
but to use this result to show that 
$\phi$ is completely positive, a detailed analysis of the case 
where $\phi$ is a homogeneous polynomial is required. Here 
the isomorphism $C^*(S^n(\cA), p_n) \cong S^n(C^*(\cA,p))$ is a key ingredient. 
All these partial results are combined in Section~\ref{sec:6} to a proof of 
Theorem~\ref{thm:intro}.

\section{Basic concepts
% and the main theorem
} \mlabel{sec:2}

In this section we introduce real seminormed involutive algebras $(\cA,p)$ and their 
enveloping $C^*$-algebras $C^*(\cA,p)$. We further introduce completely positive maps 
on these algebras and their open unit balls. 
We conclude this section with a discussion of 
complexifications and embeddings of non-unital algebras into unital ones and how 
this is reflected by the enveloping $C^*$-algebra. 

\subsection*{Seminormed involutive algebras} 

\begin{defn}\label{basic2}
A {\it real involutive algebra} is a real associative algebra $\cA$, 
endowed with an involutive antiautomorphism $a \mapsto a^*$. 

(a) A submultiplicative seminorm $p \: \cA
\to [0,\infty)$ is called {\it involutive} 
if $p(a^*) = p(a)$ for all $a \in \cA$. 
If $\cA$ is unital, we further assume $p(\1) = 1$. 
A {\it seminormed involutive algebra} is a 
pair $(\cA,p)$, consisting of an involutive algebra $\cA$ 
and an involutive submultiplicative 
seminorm $p$ on $\cA$ (cf.\ \cite[Def.~III.2.3]{Ne00}).   

(b) We will use the following notation: 
\begin{itemize}
\item $\U(\cA):=\{u\in\cA\mid u^*u=uu^*=\1\}$ for the \emph{unitary group} of a unital 
$*$-algebra $\cA$ and for a non-unital algebra we put $\U(\cA):=\U(\cA^1)\cap(\1+\cA)$, 
where $\cA^1 = \cA \oplus \R \1$ is the unitization of $\cA$; 
\item $\ball(\cA,p):=\{a\in\cA\mid p(a)< 1\}$ for the \emph{open ball semigroup} of $(\cA,p)$; 
\item $\oline{\ball}(\cA,p):=\{a\in\cA\mid p(a)\le1\}$ for the \emph{closed ball semigroup} of $(\cA,p)$.   
\end{itemize}
Both $\ball(\cA,p)$ and $\oline{\ball}(\cA,p)$ are $*$-semigroups, 
and if $\cA$ is unital, then $\1 \in\oline{\ball}(\cA,p)$. 
We will omit $p$ from this  notation whenever $p$ is clear from the context, 
for instance if $\cA$ is a $C^*$-algebra with the norm $p$. 
If $\cH$ is a complex Hilbert space and $B(\cH)$ is the von Neumann algebra of all 
bounded linear operators on $\cH$, 
then we also write $C(\cH):=\oline{\ball}(B(\cH))$. %??? where used

(c) We define the \emph{positive cone} of $\cA$ as $\cA_+:=\conv(\{a^*a\mid a\in\cA\}),$ 
where $\conv$ denotes the convex hull. 
\end{defn}

Here are some examples of seminormed involutive algebras. 
\begin{ex}\label{functions} (a) Besides $C^*$-algebras and 
Banach $*$-algebras, 
a typical example for the above setting is 
the Fr\'echet algebra of smooth $\cA_0$-valued functions $\cA=\cC^\infty(M,\cA_0)$  
on any compact manifold~$M$, where $\cA_0$ is any Banach $*$-algebra and 
$p(a) = \|a\|_\infty$. 
In this case $C^*(\cA)=\cC(M,\cA_0)$ is the $C^*$-algebra of continuous $\cA_0$-valued functions on $M$ 
and $\eta_{\cA}\colon \cC^\infty(M,\cA_0)\hookrightarrow\cC(M,\cA_0)$ is the isometric inclusion map. 

(b) Let $(S,*,\alpha)$ be an involutive semigroup $(S,*)$, endowed with an 
 {\it absolute value}, i.e., a non-negative function $\alpha \: S \to [0,\infty)$ 
satisfying
\[ \alpha(s^*) = \alpha(s) \quad \mbox{ and }\quad 
\alpha(st) \leq \alpha(s) \alpha(t) \quad \mbox{ for } \quad s,t \in S.\] 
Then the real semigroup algebra $\cA := \R[S]$ inherits an involution by 
the linear extension of the involution on $S$. 
Writing $\delta_s \in \cA$ for the basis element corresponding to $s \in S$, we obtain by 
\[ p\Big(\sum_{s \in S} c_s \delta_s\Big) 
:= \sum_{s \in S} |c_s| \alpha(s) \] 
a seminormed real involutive algebra $(\R[S],p)$.  (cf.\ \cite[Lemma~III.2.4]{Ne00}). 
\end{ex}

\begin{defn} \mlabel{def:1.1a} (Enveloping $C^*$-algebras)
For any seminormed involutive algebra $(\cA,p)$, there exist a 
$C^*$-algebra $C^*(\cA,p)$ and a morphism 
$\eta_\cA \: \cA \to C^*(\cA,p)$ of involutive algebras satisfying 
$\|\eta_\cA(a)\| \leq p(a)$ with the universal property that, 
for every real linear homomorphism $\beta \: \cA \to \cB$, where $\cB$ is a 
$C^*$-algebra and $\|\beta(a)\| \leq p(a)$ for every $a \in \cA$, there exists a 
unique morphism $\tilde\beta \: C^*(\cA,p) \to \cB$ of $C^*$-algebras satisfying 
$\tilde\beta \circ \eta_\cA = \beta$ 
(see \cite[Sect.~III.2]{Ne00}). The morphism $\eta \: \cA \to C^*(\cA,p)$ 
of real seminormed involutive algebras is called the {\it enveloping $C^*$-algebra of 
$(\cA,p)$}. 
\end{defn}

\begin{rem} \mlabel{rem:cstar} 
(a) Once we have fixed the involutive seminorm $p$, it determines a 
topology on $\cA$ so that mostly it is not necessary to specify any finer topology. 
However, for any topology~$\tau$ on $\cA$ which is finer than the 
one specified by $p$,  the 
homomorphism $\eta_\cA \: \cA \to C^*(\cA,p)$ is continuous, so that every 
linear representation of $\cA$ satisfying $\|\pi(a)\| \leq p(a)$ for $a \in\cA$, 
is $\tau$-continuous. 

(b) A linear homomorphism $\pi \: \cA \to \cB$, where $\cB$ is a $C^*$-algebra 
satisfies $\|\pi(a)\| \leq p(a)$ for every $a \in \cA$ if and only if 
$\pi(\ball(\cA,p))$ is bounded (cf.\ Remark~\ref{SzN-rem}(b) below). 
This condition implies that $\|\pi(a)\| \leq q(a)$ holds for every equivalent 
submultiplicative involutive seminorm on $\cA$. It follows in particular, that 
$C^*(\cA,p) \cong C^*(\cA,q)$ whenever there exist constants $c_1, c_2 > 0$ 
such that $c_1 q \leq p \leq c_2 q$. 

(c) Since $p$ is multiplicative on $\cA$, the subspace 
$\cI_p := \{ a \in \cA \: p(a) =0\}$ is an ideal of $\cA$, and the 
quotient $\cA/\cI_p$ inherits the structure of a real normed involutive algebra. 
Accordingly, its completion is a real Banach $*$-algebra $\cA_p$ whose 
enveloping $C^*$-algebra coincides with $C^*(\cA,p)$. 

We also point out that if $\cA$ is a $C^*$-algebra, then $C^*(\cA)$ defined as above 
is not equal to $\cA$. It is isomorphic to $\cA_{\C}\simeq\cA\oplus\overline{\cA}$, 
since the universal property from Definition~\ref{def:1.1a} involves all the real 
linear homomorphisms $\beta$, not only the complex linear ones.
\end{rem}

\begin{ex} \mlabel{ex:abel} Suppose that $\cA$ is commutative and let 
\[ \hat\cA_p := \{ \chi \in\Hom(\cA,\R) \setminus \{0\} \: \vert\chi\vert \leq p\},\] 
where $\Hom(\cA,\R)$ denotes the set of all $*$-homomorphisms $\cA \to \R$. 
Note that every such homomorphism annihilates all skew-symmetric elements, so that we 
may assume that $a^* = a$ for every $a \in \cA$. Then the weak-$*$-topology turns 
$\hat\cA_p$ into a locally compact space and it is easy to see that the natural map 
\[ \eta \: \cA \to C_0(\hat \cA_p) = C_0(\hat \cA_p,\C), 
\quad \eta(a)= \hat a, \quad \hat a(\chi) :=\chi(a) \] 
is the enveloping $C^*$-algebra of $(\cA,p)$, so that 
$C^*(\cA,p) \cong C_0(\hat \cA_p)$. 
\end{ex}

\begin{ex} (a) Consider $\cA = M_n(\K)$ for $\K \in \{\R,\C,\H\}$ with the natural involution 
$(a_{ij})^* = (\oline{a_{ji}})$ and any involutive seminorm $p$ with 
$\|a\| \leq p(a)$ for $a\in \cA$. 
Since $\cA$ is simple, the ideal $\cI_p := \{ a \in \cA \: p(a) = 0\}$ is trivial, 
and since all non-degenerate representations of $\cA$ are multiples of the complexification 
of the identical representations, we obtain 
\[ C^*(\cA,p) \cong \cA_\C \cong M_n(\K_\C) \cong 
\begin{cases}
  M_n(\C) & \text{ for } \K = \R \\ 
  M_n(\C) \oplus M_n(\C) & \text{ for } \K = \C \\ 
  M_n(M_2(\C)) \cong M_{2n}(\C) & \text{ for } \K = \H.
\end{cases}\]

(b) For $\cA := M_\infty(\K) := \bigcup_{n \in\N} M_n(\K)$, one thus obtains for any
$p\geq \|\cdot\|$ that 
\[ C^*(\cA,p) \cong \indlim C^*(M_n(\K),p) \cong K(\ell^2(\N,\K)_\C) 
\cong 
\begin{cases}
  K(\ell^2) & \text{ for } \K = \R \\ 
  K(\ell^2) \oplus K(\ell^2) & \text{ for } \K = \C \\ 
  K(\ell^2) & \text{ for } \K = \H.
\end{cases}\]
\end{ex}

\subsection*{Completely positive functions} 

\begin{defn}\mlabel{defCP} (a) A map $\phi \: \cA \to \cB$ between involutive algebras is said to be 
{\it positive} if $\phi(\cA_+) \subeq \cB_+$, i.e., positive elements of $\cA$ are mapped 
to positive elements of $\cB$. 
We call $\phi$ {\it completely positive} if, for every $n \in \N_0$, 
the induced map 
\[ \phi_n \: M_n(\cA) \to M_n(\cB), \quad 
(a_{ij}) \mapsto (\phi(a_{ij})) \] 
is positive, i.e., for $A = (a_{ij}) \in M_n(\cA)_+$, the  matrix 
$\phi_n(A)  \in M_n(\cB)$ is positive. 
This means that, for $B_1, \ldots, B_m \in  M_n(\cA)$, we have 
$\phi_n\big(\sum_{\ell=1}^m B_\ell^* B_\ell\big) \geq 0$. In view of 
\[ \Big(\sum_{\ell=1}^m B_\ell^* B_\ell\Big)_{ij}
= \sum_{\ell, k} B_{\ell, ki}^* B_{\ell,kj},\] 
this is equivalent to the requirement that 
\begin{equation}
  \label{eq:matprod}
\Big(\phi\Big(\sum_{k = 1}^m a_{ki}^* a_{kj}\Big)\Big)_{1 \leq i,j \leq n} \in M_n(\cB)_+ 
\quad \mbox{ for } \quad (a_{ki}) \in M_{m,n}(\cA), m \in \N. 
\end{equation}

(b) This concept can be extended to maps defined on $\ball(\cA,p)$ as follows. 
We call a map 
$\vphi  \: \ball(\cA,p) \to \cB$ {\it completely positive}  
if $A = (a_{ij}) \in M_n(\cA)_+$ and 
$a_{ij} \in \ball(\cA,p)$ imply that the matrix $\phi_n(A) = (\vphi(a_{ij})) 
\in M_n(\cB)$ is positive, 
for every $n\geq 1$. 
%(c) For the semigroup $C(\cA,p)$ we define complete positivity of types (W/S) accordingly. 
 \end{defn}

\begin{defn} \mlabel{def:3.1}
Let $S$ be a $*$-semigroup and $\cB$ be a $*$-algebra. 
A function $\varphi\colon S\to\cB$ is said to be {\it positive definite} 
if, for every $n \in \N$ and $s_1,\dots,s_n\in S$, 
the matrix $(\varphi(s_j^*s_k))\in M_n(\cB)$ 
is positive. 
\end{defn}

The following lemma shows that completely positive functions are particular 
positive definite functions and one of the main points of the present 
paper is to explain in which sense they are particular among the 
positive definite functions. 

\begin{lem}\mlabel{lem:3.2}
If $\phi \: \ball(\cA,p) \to \cB$ is completely positive, 
then $\phi$ is positive definite. 
\end{lem}

\begin{prf} If $a_1, \ldots, a_n \in \ball(\cA,p)$ and 
$A = (\delta_{1i} a_j)_{1 \leq i,j\leq n}$, then 
$A^* A  = (a_i^* a_j)_{1 \leq i,j \leq n}$ with 
$p(a_i^*a_j) < 1$. This implies the lemma. 
\end{prf}

\begin{rem} \mlabel{rem:1.5} (a) If $\phi \: \cA \to \cB$ and $\psi \: \cB \to \cC$ are completely 
positive maps, then so is their composition $\psi\circ \phi \: \cA \to \cC$. 
The same conclusion holds if 
$\phi \: \ball(\cA,p) \to \ball(\cB,q) \subeq \cB$ 
and $\psi \: \ball(\cB,q) \to \cC$ are completely positive, or if 
$\vphi \: \ball(\cA,p) \to \cB$ is completely positive and 
$\psi \: \cB \to \cC$ is completely positive. 

(b) If $\vphi \: \cA \to \cB$ is a linear homomorphism of $*$-algebras, 
then so are the corresponding maps $\vphi_n \: M_n(\cA) \to M_n(\cB)$, 
and this implies that $\vphi$ is completely positive. 

(c) From \eqref{eq:matprod} it follows that a 
linear map $\phi \: \cA \to \cB$ is completely positive if and only if 
it is positive definite (cf.\ \cite[Lemma~3.13]{Pau02}; \cite[p.~112]{Arv10}). 
In particular, a linear functional $\phi \: \cA \to \C$ is completely positive if and only 
if it is positive (cf.\ \cite[Thm.~3]{St55} for the case where $\cA$ is a unital $C^*$-algebra). 

(d) If $(S,*)$ is an involutive semigroup, 
then every map $\vphi \: S \to \cB$ extends to a linear map $\phi_L \: \R[S] \to \cB$. 
Using (c), it is easy to see that $\phi_L$ is positive definite if and only 
if $\phi$ is positive definite. 
 \end{rem}

The next lemma shows in particular that a homogeneous function 
on $\cA$ is completely positive if and only if it is so on 
$\ball(\cA,p)$. 

\begin{lem} \mlabel{lem:1.5} 
Suppose that the map $\vphi \: \ball(\cA,p) \to \cB$ 
is homogeneous of degree~$\alpha\in\R$, i.e.,  
for $0 < r < 1$ and $a \in \ball(\cA,p)$, we have 
$\vphi(ra)=r^\alpha \vphi(a)$. 
Then $\vphi$ extends by 
\[ \hat\vphi(ra) := r^\alpha \vphi(a) \quad \mbox{ for } \quad 
0 < r\text{ and } a \in \ball(\cA,p)\] 
to a map $\hat \vphi \: \cA \to \cB$. 
Then $\vphi$ is completely positive if and only if $\hat\vphi$ is completely positive. 
 \end{lem}

\begin{prf} We first observe that $\hat\vphi$ is well-defined.
It is clear that $\phi$ is completely positive if $\hat\phi$ is. 
To see the converse, we note that, 
for any matrix $A = (a_{ij}) \in M_n(\cA)_+$, there exists an $r > 0$ with 
$r p(a_{ij}) < 1$ for every $i,j$. Then 
$\hat \phi_n(A) = r^{-\alpha} \hat\phi_n(\alpha A) \in M_n(\cB)_+ $ 
implies that $\hat\vphi \: \cA \to \cB$ is completely positive. 
\end{prf}

\subsection*{Complexification of a real involutive algebra} 

\begin{defn} \mlabel{def:complexif} 
(Complexification of a seminormed involutive algebra) 
Let $(\cA,p)$ be a real seminormed involutive algebra. 
Then its complexification $\cA_\C$ 
is an involutive algebra in the usual sense with respect to the 
antilinear extension $(x+iy)^* := x^* - i y^*$, for  $x,y \in \cA$, of the involution. 
Further, 
\[ p_\C(a + ib) := p(a) + p(b) \] 
is a seminorm satisfying 
\[ p_\C((a+ib)^*) = p_\C(a^* - i b^*) = p(a^*) + p(b^*) = p(a) + p(b) = p_\C(a + ib) \] 
and 
\begin{align*}
p_\C((a+ib)(c+ id)) 
&= p(ac -bd) + p(bc + ad) 
\leq p(a)p(c) + p(b) p(d) + p(b) p(c) + p(a) p(d) \\
&= (p(a) + p(b))(p(c) + p(d)) = p_\C(a+ib)p_\C(c+id),
\end{align*}
so that $(\cA_\C, p_\C)$ is a seminormed involutive algebra. We call it the {\it complexification 
of $(\cA,p)$}. 
\end{defn}

\begin{lem} For a real involutive algebra $(\cA,p)$, the complex linear extension 
\[ \eta_{\cA}^\C \: \cA_\C \to C^*(\cA,p), \quad a + ib \mapsto \eta_\cA(a) + i \eta_\cA(b) \] 
has the universal property of the enveloping $C^*$-algebra in the category of complex involutive algebras, 
i.e., for every complex linear $*$-homomorphism $\beta \: \cA_\C \to \cB$ to a $C^*$-algebra 
satisfying $\|\beta(a)\| \leq p_\C(a)$ for $a \in \cA_\C$, there exists a 
unique morphism of $C^*$-algebras $\hat\beta \: C^*(\cA,p) \to \cB$ with 
$\hat\beta \circ \eta_\cA^\C = \beta$. In this sense we have 
$C^*(\cA_\C, p_\C) \cong C^*(\cA,p)$. 
\end{lem}

\begin{prf}  Let $\beta \: \cA_\C \to \cB$ 
be a complex linear $*$-homomorphism to a $C^*$-algebra 
satisfying $\|\beta(a)\| \leq p_\C(a)$ for $a \in \cA_\C$. 
Then $\|\beta(a) \| \leq p(a)$ for $a \in \cA$, so that the universal property of 
$\eta_\cA$ implies the existence of a morphism of $C^*$-algebras $\hat\beta \: C^*(\cA,p) \to \cB$ 
with $\hat\beta \circ \eta_\cA = \beta\res_\cA$. Now complex linear extension yields 
$\hat\beta \circ \eta_\cA^\C = \beta$. 
\end{prf}

\begin{rem}
 Note that a representation $\beta \: \cA \to B(\cH)$ satisfies 
$\|\beta(a)\| \leq p(a)$ for every $a \in \cA$ if and only if its complex linear extension 
$\beta^\C$ satisfies $\|\beta^\C(a + ib)\| \leq p_\C(a+ ib)$ for $a,b \in \cA$. 
\end{rem} 

\subsection*{Unitization of seminormed involutive algebras}
%{Embedding seminormed involutive algebras into unital ones} 

\begin{defn} \mlabel{def:A1}
Let $(\cA,p)$ be a real seminormed involutive algebra 
and $\cA^1 := \cA \oplus \R$, $\1 := (0,1)$, the corresponding unital algebra 
with the involution $(a,\lambda)^* := (a^*,\lambda)$. 
Then $p^1(a,\lambda) := p(a) + |\lambda|$ turns $(\cA^1, p^1)$ into a unital 
seminormed involutive algebra. 
\end{defn}

\begin{rem}\mlabel{rem:unitization}
Every $*$-homomorphism $\beta \: \cA \to \cB$ into a unital $C^*$-algebra $\cB$,  
satisfying the condition   
$\|\beta(a)\| \leq p(a)$,  extends by 
$\beta^1(a,\lambda) := \beta(a) + \lambda\1$ to a morphism of unital 
$*$-algebras satisfying $\|\beta^1(a)\| \leq p^1(a)$ for 
$a \in \cA^1$. This implies that 
$C^*(\cA^1, p^1)$ has the universal property of the {\it enveloping 
unital $C^*$-algebra} of $(\cA,p)$, i.e., 
\[ C^*(\cA,p)^1 \cong C^*(\cA^1,p^1), \] 
(cf.\ \cite[Lemma 2.10]{Bha98} for the special case of Banach algebras). 
\end{rem}

\begin{ex} \label{ex:2.17} (a) For every $c \in \cA$, the map 
$T_c \: \cA^1 \to \cA, a \mapsto c^*ac$ is completely positive. 
For $n \in \N$, the induced map $T_{c,n} \: M_n(\cA^1) \to M_n(\cA)$ 
is given by $T_{c,n} = T_C$ with $C = c E_n$, where $E_n\in M_n(\cA)$ is the identity 
matrix. 

For the sake of later reference, we note that this implies in particular that  
if $A = (a_{k\ell}) \in M_n(\cA^1)_+ \subeq M_n(M_N(\cA^1)) \cong M_{nN}(\cA^1)$ 
and $s_1, \ldots, s_N \in \cA$, then the matrix 
$(s_i^* a_{k\ell} s_j) \in M_{nN}(\cA)$ is positive because 
it can be written as $S^* A S$ for 
$S = (\delta_{ij} s_j)_{1 \leq i,j \leq N, 1 \leq k,\ell \leq n}$.

(b) Fix $s_1, \ldots, s_N \in \cA$. We claim that the map 
\[ \beta \: M_N(\cA^1) \to \cA, \quad 
\beta((a_{ij}) := \sum_{i,j=1}^N s_i^* a_{ij} s_j.\] 
is completely positive. The corresponding map 
\begin{eqnarray}\label{eq:betan}
\beta_n \: M_n(M_N(\cA^1)) \cong M_{nN}(\cA^1)  &\to& M_n(\cA), \notag \\
\beta_n((a_{(i,k),(j,\ell)})) &
= &\Big(\sum_{i,j=1}^N s_i^* a_{(i,k),(j,\ell)} s_j\Big)_{1 \leq k,\ell \leq n}\
\end{eqnarray}
is of the same structure with $s_i$ replaced by 
the matrix $s_i E_n \in M_n(\cA)$. Therefore it suffices to show that 
$\beta$ is positive. As $\beta$ is linear, this means that 
$\beta(B^*B) \geq 0$ for $B \in M_N(\cA_1)$. 
For $A = B^*B$, we have 
\[ \beta(B^*B) 
= \sum_{i,k,j=1}^N s_i^* b_{ki}^* b_{kj} s_j 
= \sum_{k = 1}^N \Big(\sum_{j=1}^N b_{kj} s_j\Big)^*\Big(\sum_{j=1}^N b_{kj} s_j\Big) \geq 0.\]
\end{ex}

\section{Positive definite functions and dilations} 
\mlabel{sec:3}

We start this section by a discussion of the representation 
$(\pi_\phi, \cH_\phi)$ associated to a (bounded) positive definite 
function $\phi \: S \to B(V)$ on a $*$-semigroup $S$. 
Here $\cH_\phi$ is a subspace of the space 
$V^S$ of $V$-valued functions on $S$, on which $S$ acts by right translations. 
We then discuss dilations of positive definite functions and criteria for their 
existence, such as the existence of a unit in $S$. 
For non-unital 
semigroups, the existence of dilations 
%becomes 
is a subtle issue. 
Based on the concept of a structured $*$-semigroup and multiplier techniques, we 
prove in Theorem~\ref{thm:Arv2.2} the existence of dilations for positive definite 
functions on the canonical ideal $S^\circ$ of a structured $*$-semigroup, 
and this applies in particular to the subsemigroup $S^\circ = \ball(\cA,p)$ of 
$S = \oline\ball(\cA,p)$ when $\cA$ is unital. 
Generalizing the classical Stinespring Dilation Theorem (\cite{St55}), we 
prove in Proposition~\ref{cp-univ} that a linear completely positive function 
$\phi \: \cA \to B(V)$ which is bounded on $\ball(\cA,p)$ factors through 
$\eta_\cA \: \cA \to C^*(\cA,p)$ if and only if it has a dilation.

\subsection{Preliminaries on positive definite functions}

Throughout this subsection, $(S,*)$ is an involutive semigroup. 

\begin{rem}\label{nonsense} (The GNS representation associated to $\phi$) 
We recall from \cite{Ne00} 
some basic facts on the reproducing kernel Hilbert space 
associated with any positive definite function 
$\varphi \colon S \to B(V)$, where $S$ is any $*$-semigroup. 
Define the positive definite kernel 
\[ K=K_\varphi\colon S\times S\to B(V)\quad \mbox{ by } \quad K(s,t) := \varphi(st^*).\]  
The corresponding reproducing kernel space $\cH_\vphi \subeq V^S$ 
($V^S$ stands for the vector space of all $V$-valued functions on~$S$), 
can be constructed as follows. 
For $t \in S, v \in V$, we put $K_{t,v} := K(\cdot, t)v \in V^S$. Then 
\[ \cH_\varphi^0:=\Spann\{K_{t,v}  \mid t\in S,\ v\in V\}  \] 
is a pre-Hilbert space with respect to a  scalar product specified by 
\[ \la K_{s,v}, K_{t,w} \ra = \la K(t,s)v,w\ra_V =\la\varphi(ts^*)v,w\ra_V 
\quad \mbox{ for } \quad 
s,t \in S,v,w \in V. \]
The completion $\cH_\vphi$ of $\cH_\vphi^0$ 
has a natural inclusion into $V^S$ with continuous evaluation operators 
\[ K_s \: \cH_\vphi \to V, f \mapsto f(s) \qquad \mbox{ satisfying} \quad 
K_s K_t^* = K(s,t) \quad \mbox{ for } \quad s,t \in S.\] 
In particular, we have 
$K_{t,v} = K_t^* v$ for $t \in S, v \in V$. 

As in \cite[Sect. II.3]{Ne00}, we denote by $B^0(\cH_\varphi^0)$ 
the $*$-algebra consisting of the linear operators $A\colon\cH_\varphi^0\to\cH_\varphi^0$ 
for which there exists a linear operator $A^\sharp\colon\cH_\varphi^0\to\cH_\varphi^0$ with 
$\la Av,w\ra=\la v,A^\sharp w\ra$ for all $v,w\in\cH_\varphi^0$. 
Then we obtain a $*$-representation 
$\pi_\varphi^0 \: S \to B^0(\cH_\varphi^0)$ given by  
$(\pi_\varphi^0(s)f)(t) := f(ts)$ for all $t,s \in S$. 
For the operator closures 
$\pi_\phi(s) := \oline{\pi_\phi^0(s)}$, we then obtain 
\begin{equation}
  \label{eq:k-rels}
 K_t \pi_\vphi(s) = K_{ts} \quad \mbox{ and } \quad \pi_\vphi(s) K_t^* = K_{ts^*}^*
\quad 
\mbox{ for } \quad t,s \in S.
\end{equation}
\end{rem}

\begin{rem} If $\vphi \: S \to B(V)$ is constant, then all functions 
$K_{t,v}(s) = K(s,t)v = \vphi(st^*)v$ are constant, so that the representation 
$(\pi_\vphi, \cH_\vphi)$ is trivial. 

If, conversely, the representation $\pi_\vphi$ is trivial, 
then $K_t = K_{ts}$ for $t,s \in S$ by \eqref{eq:k-rels}. 
If $S$ has a unit, this implies that $\vphi$ is constant. 
\end{rem}

\begin{defn}\label{bddness}
Let $S$ be a $*$-semigroup, $V$ be a complex Hilbert space 
and $\phi \: S \to B(V)$ be a positive definite function. 

(a) If $\alpha \: S \to [0,\infty)$ is an absolute value 
(cf.\ Examples~\ref{functions}(b)), then $\phi$ is said to be 
{\it $\alpha$-bounded} if 
\[ \varphi(s^*t^*ts)\le \alpha(t)^2\varphi(s^*s) \quad \mbox{ for } \quad 
t,s\in S.\] 
This condition is equivalent to $\|\pi_\phi(t)\| \leq \alpha(t)$ for the corresponding 
GNS representation $(\pi_\phi, \cH_\phi)$ (cf.\ \cite[Th.~III.1.3]{Ne00}; 
and \cite{Seb86} for related conditions). 

(b)  A triple $(\pi,\cH,\iota)$ 
consisting of a $*$-representation $\pi\colon S\to B(\cH)$ on 
the complex Hilbert space $\cH$ and a bounded linear operator 
$\iota\colon V\to\cH$ (called the {\it linking operator}), 
is called a {\it dilation of $\phi$} 
if 
\[ \phi(s) = \iota^* \pi(s) \iota \quad \mbox{ for } \quad s \in S.\] 
A dilation is said to be {\it minimal} if 
$\pi(S)\iota(V)$ spans a dense subspace of $\cH$.
A positive definite function for which a dilation exists is called 
{\it dilatable}.
\begin{footnote}  {The above dilatable functions were called Stinespring representable in \cite{Bha98} in the special case when $S=(\cA,\cdot)$ for some Banach $*$-algebra~$\cA$.}
\end{footnote}
\end{defn}

The following proposition 
shows in particular that the boundedness characterization 
of \cite[Th.~III.1.19(4)]{Ne00} can also be used for operator-valued functions. 

\begin{prop}   \mlabel{prop:minuni} 
For any positive definite function $\phi \: S \to B(V)$ on the involutive semigroup $(S,*)$, 
the following assertions hold: 
\begin{itemize}
\item[\rm(a)] If $S$ is unital and $\phi$ is $\alpha$-bounded for some 
absolute value $\alpha$ on $S$, then $(\pi_\phi, \cH_\phi, K_\1^*)$ 
is a minimal dilation of $\phi$ satisfying $\|\pi_\phi(s)\| \leq \alpha(s)$ for $s \in S$. 
\item[\rm(b)] If $(\pi, \cH, \iota)$ is a minimal dilation 
of $\phi$, then there exists a unique unitary operator 
$\Phi \: \cH \to \cH_\phi$ intertwining  $\pi$ with $\pi_\phi$ and satisfying 
$\Phi(\iota(v)) = \phi(\cdot) \cdot v$ for all $v\in V$. 
\item[\rm(c)] {\rm(Uniqueness of minimal dilations)} 
If $(\pi_j, \cH_j, \iota_j)$, $j =1,2$, are two minimal dilations 
of the positive definite function $\phi \: S \to B(V)$, 
then there exists a unique unitary operator $U\colon\cH_1\to\cH_2$ 
intertwining $\pi_1$ with $\pi_2$ and satisfying $U \circ \iota_1 = \iota_2$.  
\end{itemize}
\end{prop}

\begin{prf} (a) We will use the notation from Remark~\ref{nonsense}. 
In view of \cite[Th.~III.1.3]{Ne00},  the $\alpha$-boundedness of $\phi$ implies 
that $\Vert\pi_\varphi^0(s)\Vert\le\alpha(s)$ for every $s\in S$.  
Since $\cH_\varphi^0$ is dense in $\cH_\varphi$, it follows that the linear operator 
$\pi_\varphi^0(s)\in B(\cH_\varphi^0)$ has a unique extension $\pi_\varphi(s)\in B(\cH_\varphi)$ 
with $\Vert\pi_\varphi(s)\Vert\le\alpha(s)$. 
This is also given by right translations 
\[ (\pi_\varphi(s)f)(t) := f(ts) \quad \mbox{ for } \quad t,s \in S, f \in \cH_\phi \]  
(cf.\ \cite[Prop.~II.4.9]{Ne00}). 
For $s \in S$ and $\iota := K_\1^*$, we have 
\[ \iota^* \pi_\phi(s) \iota(v) 
= K_\1 \pi_\phi(s) K_\1^* v = (K_\1^* v)(s) = \phi(s)v,\] 
so that $(\pi_\phi, \cH_\phi, \iota)$ is a dilation of $\phi$. 
Its minimality follows from the fact that, for 
$f \in \cH_\phi$, the vanishing of 
$f(s) = K_\1 \pi(s) f$ for every $s \in S$ implies $f = 0$. 

(b) We consider the map 
\[ \Phi \: \cH \to V^S, \quad \Phi(w)(s) := \iota^* \pi(s)w.\] 
This map is $S$-equivariant with respect to the right translation action of $S$ on $V^S$ and 
satisfies $\Phi(\iota(v)) = \phi \cdot v$ for all $v\in V$. 
If $\Phi(w) = 0$, then $\pi(S)w \subeq  \ker \iota^* = \iota(V)^\bot$ leads 
to $w \in (\pi(S)\iota(V))^\bot = \{0\}$, so that $\Phi$ is injective. 
Therefore $\Phi(\cH) \subeq V^S$ is a reproducing kernel Hilbert space whose kernel 
is given by 
\[ K(s,t) = (\iota^* \pi(s))(\iota^* \pi(t))^* = \iota^* \pi(st^*) \iota  
= \phi(st^*).\] 
We conclude that $(\pi, \cH)$ is equivalent to $(\pi_\phi, \cH_\phi)$. 

(c) follows from (b) with $U := \Phi_2^{-1} \circ \Phi_1$. 
\end{prf}

\begin{rem}\label{SzN-rem}
Let $S$ be a $*$-semigroup, $\varphi\colon S\to B(V)$ be a positive definite function 
and $\alpha \: S \to [0,\infty)$ be an absolute value. 

(a) If $\|\varphi(s)\| \leq M \alpha(s)$ for all $s\in S$ and some $M \geq 0$, 
then $\phi$ is $\alpha$-bounded (see the supplement to \cite[Appendix 1]{RSzN72} 
or \cite[Cor.~III.1.20]{Ne00}).
In particular, Proposition~\ref{prop:minuni}(a)  applies to 
every bounded positive definite function $\varphi\colon S\to B(V)$ 
with the constant absolute value $\alpha = 1$.
\begin{footnote}
{In the case $V=\C$ one thus recovers \cite[Lemma 2.3]{Arv87}.}
\end{footnote}

(b) If $\pi\colon S\to B(\cH)$ is a $*$-representation 
for which $\sup\limits_{s\in S}\Vert\pi(s)\Vert<\infty$, then actually 
$\sup\limits_{s\in S}\Vert\pi(s)\Vert\le 1$ because, for $s \in S$, 
the boundedness of the sequence 
\[ \|\pi(s)\|^{2n} =  \|\pi(s^*s)\|^{n} =  \|\pi((s^*s)^n)\| \] 
is equivalent to $\|\pi(s)\|  \leq 1$. 
\end{rem}

\begin{lem}\label{Arv2.2-cor2}
Let $S$ be a $*$-semigroup with a zero element 
$0$, i.e., $0 \cdot s = s \cdot 0 = 0$ for $s \in S$. 
If $\varphi\colon S\to B(V)$ is an $\alpha$-bounded positive definite function, 
then the function $\varphi-\varphi(0)$ is also positive definite. 
\end{lem}

\begin{prf} Since $0$ is a zero element, the same holds for $0^*$, and this leads to 
$0 = 0^* 0 = 0^*$. Therefore the operator $P:=\pi_\phi(0)$ is a hermitian projection 
satisfying $P\pi_\phi(s)=\pi_\phi(s)P=P$ for every $s\in S$. It acts on 
$\cH_\phi$ by 
\[ (\pi_\phi(0)f)(s) = f(s0) = f(0).\] 
This implies that $\cH_0 := P(\cH_\phi)$ is the subspace of constant functions 
contained in $\cH_\phi$, and that $\cH_0^\bot = \{ f \in \cH_\phi \: f(0) = 0\}$. 
Accordingly, the reproducing kernel $K(s,t) = \phi(st^*)$ decomposes as 
\[ K_s = \ev_0 \oplus (K_s - \ev_0) \in B(\cH_0, V) \oplus B(\cH_0^\bot,V) 
\cong B(\cH_\phi,V). \] 
This leads to 
\begin{align*}
(K_s - \ev_0)(K_t - \ev_0)^*
&= K_s K_t^* - \ev_0 K_t^* - K_s \ev_0^* + \ev_0 \ev_0^* \\
&= \phi(st^*) - \phi(0) - \phi(0) + \phi(0) = \phi(st^*) - \phi(0).
\end{align*}
We conclude that $\phi - \phi(0)$ is positive definite. 
\end{prf}

\begin{prop} \mlabel{prop:exis-dil} {\rm(Existence of minimal dilations)} 
Consider the unital 
involutive semigroup $S^1 := S \dot\cup\{\1\}$ with $\1^* = \1$ 
and an absolute value $\alpha$ on $S^1$. 
For an $\alpha$-bounded 
positive definite function $\phi \: S \to B(V)$, the following are equivalent: 
\begin{itemize}
\item[\rm(i)] $\phi$ is dilatable. 
\item[\rm(ii)] For every $v \in V$, the function $\phi \cdot v \:  S \to V$ 
belongs to $\cH_\phi$. 
\item[\rm(iii)] $\phi$ extends to an 
$\alpha$-bounded positive definite function on $S^1$. 
\end{itemize}
\end{prop}

\begin{prf} (i) $\Rarrow$ (ii) follows from Proposition~\ref{prop:minuni}(b). 

(ii) $\Rarrow$ (i): Suppose, 
conversely, that $\phi \cdot v \in \cH_\phi$ for every $v \in V$. 
Then we obtain a linear map $\iota \: V \to \cH_\phi, v \mapsto \phi \cdot v$. 
Since all compositions $\ev_s \circ \iota \: V \to V, v \mapsto \phi(s)v$, 
are continuous,
the Closed Graph Theorem implies that $\iota$ is continuous. 
We further have, for $v, w \in V$, the relation  
\[ \la \iota^* \pi_\phi(s^*) \iota(v), w \ra 
= \la \pi_\phi(s^*) \phi \cdot v, \phi \cdot w \ra 
= \la K_{s}^* v, \phi \cdot w \ra 
= \la v, \phi(s) w \ra.\] 
This implies that $\phi(s) =  (\iota^* \pi_\phi(s^*) \iota)^* = \iota^* \pi_\phi(s) \iota$,
so that $(\pi_\phi, \cH_\phi, \iota)$ is a dilation of $\phi$. 

(i) $\Rarrow$ (iii): If $(\pi, \cH, \iota)$ is a dilation of $\phi$ and 
$\pi^1(\1) := \1$ is the canonical extension 
of $\pi$ to a representation of $S^1$, then 
$\phi^1(s) := \iota^* \pi^1(s) \iota$ 
is an extension of $\phi$ to an $\alpha$-bounded positive definite 
function on $S^1$. 

(iii) $\Rarrow$ (i): If an $\alpha$-bounded positive definite extension $\phi^1$ to 
$S^1$ exists, then Proposition~\ref{prop:minuni}(a) implies that 
$\phi^1$, and hence also $\phi$, is dilatable. 
\end{prf}

\begin{ex} \mlabel{ex:3.8} 
(a) If $\phi \: S \to B(V)$ is a $*$-representation, then 
$(\phi, V,\id_V)$ is a dilation of $\phi$. 

(b) If $\phi \: S \to B(V)$ is positive definite and 
$s \in S$, then $\phi_s(t) := \phi(sts^*)$ defines a positive definite 
function on $S$ which extends by $\phi_s(\1) := \phi(ss^*)$ to a 
positive definite function on $S^\1$. Therefore $\phi_s$ is dilatable
by Proposition~\ref{prop:exis-dil}. More explicitly, the 
representation $(\pi_\phi, \cH_\phi)$ satisfies 
\[  K_{s} \pi_\phi(t) K_{s}^* = K_{st} K_s^* = \phi(sts^*) = \phi_s(t),\] 
so  that $(\pi_\phi, \cH_\phi, K_s^*)$ is a dilation of $\phi_s$. 

(c) Suppose that $S$ is endowed with a topology for which the bounded positive definite 
function $\phi \: S \to B(V)$ is weak-operator continuous and for which 
there exists a right approximate identity $(\delta_j)_{j \in J}$ in $S$, i.e., 
$\lim_j s\delta_j = s$ for every $s \in S$.  

For $v \in V$, we then have 
$\| K_{\delta_j,v}\|^2 = \la \phi(\delta_j\delta_j^*)v,v\ra  \leq \|\phi\|_\infty \|v\|^2$, so that the net $(K_{\delta_j,v})_{j \in J}$ has a weak cluster point 
$\psi \in \cH_\phi$. Let $K_{\delta_{j_k},v}$ be a subnet converging to $\psi$. 
Then we obtain for $w \in V$ the relation
\[ \la \psi(s), w\ra 
= \lim_k \la K_{\delta_{j_k,v}}, K_{s,w} \ra 
= \lim_k \la \phi(s \delta_{j_k})v,w \ra = \la \phi(s)v,w\ra,\] 
so that $\psi = \phi \cdot v \in \cH_\phi$. Therefore $\phi$ is dilatable by 
Proposition~\ref{prop:exis-dil}(iii). 
\end{ex}

\begin{prop}\mlabel{cp-univ}
Let $(\cA,p)$ be a seminormed $*$-algebra and $V$ be a complex Hilbert space. 
Then $\Phi\mapsto\Phi\circ\eta_\cA$ 
is a bijective correspondence between the set of  linear completely positive 
$B(V)$-valued  maps on $C^*(\cA,p)$ 
and the dilatable completely positive functions $\phi \: \cA \to B(V)$ 
which are bounded on $\ball(\cA,p)$. 
\end{prop}

\begin{prf} 
Since $\eta_\cA$ is a linear homomorphism of $*$-algebras, 
it is completely positive (Remark~\ref{rem:1.5}) 
and its range spans a dense subspace of $C^*(\cA,p)$. 
Hence the correspondence $\Phi \mapsto \Phi \circ \eta_\cA$ 
is well defined and injective. Every completely positive map 
$\Phi \: C^*(\cA,p) \to B(V)$ extends to a 
completely positive map $\Phi^1 \: C^*(\cA,p)^1 \to B(V)$ 
by Arveson's Extension Theorem (\cite[Thm.~1.2.3]{Arv69}), 
and $C^*(\cA,p)^1 \cong C^*(\cA^1,p^1)$ by 
Remark~\ref{rem:unitization}. Therefore $\Phi$ is dilatable, 
and this implies that $\Phi \circ \eta_\cA$ is dilatable. 

Assume, conversely, that $\varphi\colon\cA\to B(V)$ is 
linear, bounded on $\ball(\cA,p)$, 
and that $(\pi, \cH, \iota)$ is a minimal dilation of $\phi$. 
Then the linearity of $\phi$ implies that the 
representation $\pi$ of $(\cA,\cdot)$ is linear. Further, 
$\|\pi(a)\| \leq p(a)$ follows from the uniqueness of minimal dilations 
(Proposition~\ref{prop:minuni}). 
By the universal property of $C^*(\cA,p)$,  there exists a unique 
$*$-homomorphism $\hat\pi\colon C^*(\cA,p)\to B(\cH)$ with $\hat\pi\circ\eta_\cA=\pi$. 
Then $\Phi(\cdot):=\iota^*\hat\pi(\cdot)\iota\colon C^*(\cA,p)\to B(V)$ is 
linear, completely positive and satisfies 
$\Phi\circ\eta_\cA=\varphi$. This completes the proof. 
\end{prf}

A Banach $*$-algebra version of Proposition~\ref{cp-univ} can be found in \cite[Th.~2.1(6)]{Bha98}. 
In that case the $*$-representations and linear positive definite $B(V)$-valued functions 
are always continuous 
(\cite[Prop.~1.3.7]{Dix64}, \cite[Cor.~A.2]{Arv10}). 
This automatic continuity property remains true 
for unital Mackey complete $*$-algebras with continuous inversion 
(\cite[Prop. 6.6]{Bi10}).

\begin{ex} \mlabel{ex:non-dil} (A non-dilatable completely positive functional) 
We consider the non-unital Banach $*$ algebra $\cA := \ell^1(\N,\R)$, 
endowed with the pointwise multiplication and $p(a) = \|a\|_1$. 
Then the enveloping $C^*$-algebra is $C^*(\cA,p) \cong c_0(\N,\C)$. 

The continuous linear functional 
\[ \phi \: \cA \to \R, \quad \phi(a) := \sum_{n = 1}^\infty a_n \] 
is a sum of characters of $\cA$, hence completely positive. 
Moreover, for $a,b \in \cA$, we have 
\[ \phi(a^*ba) \leq \|b\|_\infty \phi(a^*a) = \|b\|_\infty \|a\|_2^2,\] 
so that the representation $\pi_\phi$ of $\cA$ on 
$\cH_\phi \subeq \cA^*= \ell^\infty(\N,\C)$ extends to a representation of $C^*(\cA,p) \cong c_0(\N,\C)$ 
(Proposition~\ref{prop:minuni}). 

Note that $\phi(ab^*) = \sum_n a_n \oline{b_n}$ implies that 
$\cH_\phi \cong \ell^2(\N,\C)$, considered as a subspace of 
$\cA_\C' \cong \ell^\infty(\N,\C)$. The function 
$\phi$ itself corresponds to the constant function $1$, which 
is not contained in $\cH_\phi$. 
This implies that there exists 
no element $v \in \cH_\phi$ with $\phi(a) = \la \pi(a)v,v\ra$ for $a\in \cA$ 
(cf.\ \cite[Cor.~III.1.25]{Ne00}). We conclude that $\phi$ is 
not dilatable (Proposition~\ref{prop:exis-dil}). 

Note that $\cA$ has no bounded approximate identity because, 
in view of Example~\ref{ex:3.8}(c), this would contradict 
Proposition~\ref{prop:exis-dil}, applied to $\ball(\cA,p)$ 
(see also \cite[Rem.~IV.1.23]{Ne00}).  
\end{ex}

The preceding example, applied to $\ball(\cA,p)$, provides a 
counterexample to \cite[Th.~2]{Seb86}.

\subsection{Multipliers and structured $*$-semigroups} 

In this subsection we recall the multiplier semigroup of a $*$-semigroup 
and introduce the concept of a structured (unital) $*$-semigroup 
$(S,\gamma)$, which consists of a $*$-semigroup $S$ with a unital $*$-homomorphism 
$\gamma \: (0,1] \to S$ with central range. 
The main result of this subsection is Theorem~\ref{thm:Arv2.2} which asserts 
that bounded positive definite functions on the semigroup ideal 
$S^\circ :=\gamma((0,1))S$ are always dilatable. 

\subsubsection*{Multipliers of involutive semigroups} 

\begin{defn} (Multiplier semigroup) 
Let $(S,*)$ be an involutive
semigroup. A {\it multiplier} of $S$ is a pair $(\lambda, \rho)$  of
mappings $\lambda, \rho \: S \to S$ satisfying the following
conditions:
\begin{itemize}
\item[\rm(M1)] $a \lambda(b) = \rho(a)b$, 
\item[\rm(M2)] $\lambda(ab) = \lambda(a)b$, and 
\item[\rm(M3)] $\rho(ab) = a \rho(b).$
\end{itemize}
The map $\lambda$ is called the {\it left action} of the multiplier
and $\rho$ is called the {\it right action} of the multiplier. 
For $m = (\lambda, \rho)$ and $s \in S$ we write 
\[ ms := \lambda(s) \quad \mbox{ and } \quad sm := \rho(s).\] 
We write $M(S)$ for the set of all multipliers of $S$ and turn
it into an involutive semigroup by 
\[ (\lambda, \rho) (\lambda', \rho') := (\lambda \circ \lambda', \rho'
\circ \rho)  \quad \mbox{ and } \quad 
(\lambda, \rho)^* := (\rho^*, \lambda^*) \]
with $\lambda^*(a) := \lambda(a^*)^*$ and 
$\rho^*(a) = \rho(a^*)^*$ 
(see \cite{Jo64}).   
\end{defn}

\begin{rem} (a) 
The assignment $\eta \: S \to M(S), a \mapsto
(\lambda_a, \rho_a)$ defines a morphism of involutive semigroups which
is surjective if $S$ has an identity: in this case 
(M1) implies that $s := \lambda(\1) = \rho(\1)$, (M2) implies 
$\lambda = \lambda_s$ and (M3) that $\rho = \rho_s$. 
The multiplier semigroup acts from the right on $S$ by 
$(s,m) \mapsto sm$, 
and via $\eta$ this extends the natural right action of $S$ on itself.
This right action of $M(S)$ is related to the fact that 
$\eta(S)$ is an involutive semigroup ideal in $M(S)$. In fact, 
$\eta(s)^* = (\rho_s^*, \lambda_s^*) = (\lambda_{s^*}, \rho_{s^*}) 
= \eta(s^*),$
and
\[  (\lambda, \rho) \eta(s) 
= (\lambda \circ \lambda_s, \rho_s \circ \rho) 
= (\lambda_{\lambda(s)}, \rho_{\lambda(s)}) =
\eta\big(\lambda(s)\big). \]

(b) Let $\phi$ be a positive definite function on $S$ and $K_\phi(a,b) :=
\phi(ab^*)$ the corresponding positive definite kernel. We claim that
$K_\phi$ is invariant with respect to the action of $M(S)$ on
$S$, as is easily seen from  
\begin{align*}
K_\phi(a, bm) 
&= \phi(a(bm)^*)= \phi(am^*b^*) = K_\phi(am^*,b). 
\end{align*}

This shows that, for every positive definite function $\phi$ on $S$,
there exists a hermitian representation $\tilde\pi_\phi$ of $M(S)$ on $\cH_\phi^0$ 
satisfying $\tilde \pi_\phi \circ \eta = \pi_\phi$. 
In view of \cite[Prop.~II.2.11(ii)]{Ne00},  every 
bounded representation $(\pi, \cH)$ of $S$ yields a 
hermitian representation $\tilde\pi$ of the multiplier 
semigroup with $\tilde\pi(m) \pi(s) 
= \pi(ms)$ and therefore also 
$\tilde\pi \circ \eta = \pi$. 
To decide whether this extension acts by bounded 
operators, one needs more information on the multipliers or has to use
topological arguments involving approximate identities 
(cf.\ \cite[p.778]{FD88}). 
\end{rem}

\begin{lem} \mlabel{lem:3.13} For any non-degenerate bounded representation 
$(\pi, \cH)$ of the involutive semigroup $(S,*)$, there exists a unique 
bounded representation $(\tilde\pi,\cH)$ of its multiplier semigroup 
$M(S)$ satisfying $\tilde\pi(m) \pi(s) = \pi(ms)$ for $m \in M(S), s \in S$. 
\end{lem}

\begin{prf} Let $(\pi, \cH)$ be a non-degenerate bounded representation of $S$ 
and note that this implies that $\|\pi(s)\| \leq 1$ for every $s \in S$. 
We claim that the corresponding hermitian representation 
$\tilde\pi$ of $M(S)$ is also bounded. To see this, we may w.l.o.g.\ 
assume that $\pi$ is cyclic. So let $v \in \cH$ be an element for 
which $\pi(S)v$ is total in $\cH$. 
Then $\tilde\pi(M(S)) \pi(S)v \subeq \pi(S) v$ shows that this total 
subset is invariant under the hermitian representation 
of $M(S)$ on $\Spann(\pi(S)v)$. From $\|\pi(s)v\| \leq \|v\|$ 
it now follows that, for every $s \in S$, the positive definite 
function 
\[ \phi(m) 
:= \la \tilde\pi(m) \pi(s)v,\pi(s)v \ra 
= \la \pi(s^*ms) v, v \ra \] 
is bounded, so that \cite[Cor.~III.1.20]{Ne00} implies that 
\[ \|\tilde\pi(m)\pi(s)v\| \leq \|\pi(s)v\|.\] 
Therefore \cite[Lemma~II.3.8(iii)]{Ne00} yields 
$\|\tilde\pi(m)\| \leq 1$ for every multiplier 
$m \in M(S)$. 
\end{prf}

\begin{ex} \mlabel{ex:3.14} Let $(\cA,p)$ be a real seminormed involutive algebra and 
$(\cA^1,p^1)$ be the corresponding unital seminormed algebra (Definition~\ref{def:A1}). 
Then $\cA$ is an ideal in $\cA^1$ and 
$\ball(\cA,p)$ is a semigroup ideal in $\oline\ball(\cA^1,p^1)$, which leads to a natural homomorphism 
\[ \mu \: \oline\ball(\cA^1,p^1) \to M(\ball(\cA,p)), \quad 
\mu(m) := (\lambda_m, \rho_m).\] 
Accordingly, Lemma~\ref{lem:3.13} implies that every non-degenerate bounded 
representation $\pi$ of $\ball(\cA,p)$ extends uniquely to a representation 
$\pi^1$ of $\oline\ball(\cA^1,p^1)$ satisfying 
$\pi^1(m) \pi(a) = \pi(ma)$ for $a \in \ball(\cA,p)$, $m \in \oline\ball(\cA^1,p^1)$. 
\end{ex}

\subsubsection*{Structured $*$-semigroups}

\begin{defn}\label{struct}
A \emph{structured $*$-semigroup} is a unital $*$-semigroup $S$ 
with a homomorphism of $*$-semigroups $\gamma\colon(0,1]\to S$ 
for which $\gamma(1)=\1\in S$ and the image of $\gamma$ is contained in the center of~$S$. 
We will say that $\gamma$ is the \emph{structure homomorphism} of~$S$. 
Note that $\gamma(r)^*=\gamma(r)$ for every $r\in(0,1]$. 
We will denote 
$$(\forall r\in(0,1])(\forall x\in S)\quad xr=rx:=\gamma(r)x$$
and $S^\circ:=\bigcup\limits_{0<r<1}rS$, which is easily seen to be a $*$-subsemigroup of~$S$, 
however we may have $\1\not\in S^\circ$. 
Actually, $S^\circ$ is a semigroup ideal of $S$, 
that is, $SS^\circ\cup S^\circ S\subseteq S^\circ$. 
\end{defn}

\begin{rem} \mlabel{rem:3.16} 
If $S$ is a structured $*$-semigroup, then 
$S^\circ = S^\circ S^\circ$, as a trivial consequence of $S\cdot S = S$ ($S$ is unital). 
\end{rem}

Every $*$-semigroup with unit element is a structured $*$-semigroup 
with respect to the trivial structure homomorphism which is identically equal to the unit element. 
Here are also some less trivial examples. 

\begin{ex}\label{struct-ex0}
(a) Let $\cA$ be any unital associative $*$-algebra, 
and pick any convex self-adjoint multiplicative subset $S\subseteq\cA$ which contains $0,\1\in\cA$. 
That is, we require that $S$ should be convex, closed under the involution, 
and moreover $S\cdot S\subseteq S$ and $0,\1\in S$. 
Then the multiplicative $*$-semigroup $S$ is a typical example of a structured $*$-semigroup 
with the structure homomorphism defined by the scalar multiplication 
$\gamma\colon(0,1]\to S$, $\gamma(r)=r\1$. 

(b) %\label{struct-ex1}
Let $(\cA,p)$ is any seminormed unital $*$-algebra. 
Then the $*$-semigroup 
\[ S:=\oline{\ball}(\cA,p) := \{ a \in \cA \: p(a) \leq 1\} \] 
satisfies the requirements of (a). 
It is easily checked that 
\[ S^\circ=\ball(\cA,p) := \{ a \in \cA \: p(a) < 1\}.\] 

(c) %\label{struct-ex}
Let $\cA$ be a unital associative Banach $*$-algebra with isometric involution and 
$\|\1\| =1$. 
Then the contraction $*$-semigroup $S=\oline{\ball}(\cA) = \{ a \in \cA \: \|a\| \leq 1\}$ 
is a special case of (b) 
and it is clear that $S^\circ=\ball(\cA)$ is the strict contraction semigroup of~$\cA$. 
\end{ex}

We will prove the following generalization of \cite[Th. 2.2]{Arv87}. 
Note that it does not follows from Proposition~\ref{prop:minuni} 
because $S^\circ$ need not be unital. 

\begin{thm}\mlabel{thm:Arv2.2} {\rm(Dilation Theorem for structured  
$*$-semigroups)} 
Let $S$ be a structured $*$-semigroup and $V$ be a complex Hilbert space. 
Then every bounded positive definite function $\varphi\colon S^\circ\to B(V)$ has a dilation 
$(\pi, \cH, \iota)$ with 
\[ \|\iota\|^2 \leq \limsup_{0 < r < 1} \|\phi(\gamma(r))\|.\] 
\end{thm}

Before proving Theorem~\ref{thm:Arv2.2}, we describe some of its consequences. 

\begin{rem} \label{Arv2.2-cor1} % check! 
(a) If $(\cA,p)$ is a unital seminormed $*$-algebra, then 
$S = \oline{\ball}(\cA,p)$ is a structured $*$-semigroup 
with $S^\circ = \ball(\cA,p)$ (Example~\ref{struct-ex0}(c)). Accordingly, 
Theorem~\ref{thm:Arv2.2} applies to~$S$. 

(b) Example~\ref{ex:non-dil} shows that Theorem~\ref{thm:Arv2.2} does not 
generalize to non-unital semigroups. 
\end{rem}

For later use in the proof of Proposition~\ref{P1.2}, we 
also record the following observation. 

\begin{lem}\label{L1.1} 
If $S$ is a $*$-semigroup satisfying $S = S \cdot S$, 
 and $\varphi\colon S\to B(V)$ is an $\alpha$-bounded 
positive definite function, then the representation 
$(\pi_\varphi, \cH_\varphi)$ is non-degenerate.   
\end{lem}

\begin{prf} Let $f \in \cH_\varphi$ with $\pi_\varphi(S)f = \{0\}$. 
Then $0 =f(S \cdot S) = f(S)$ implies $f = 0$. 
\end{prf}

We now turn to the proof of Theorem~\ref{thm:Arv2.2},   
which needs the following method of extending representations 
of semigroup ideals. 
In the next proposition we use the notation $S^\circ$ from structured $*$-semigroups only in order to indicate 
how this method will be used, however this applies to semigroup $*$-ideals, 
which are more general than the ones defined by a structure homomorphism as in Definition~\ref{struct}.

\begin{prop}\label{P1.2} 
Let $S$ be a $*$-semigroup and suppose that $S^\circ \trile S$ is a semigroup $*$-ideal, i.e., 
$S S^\circ \cup S^\circ S \subeq S^\circ$, and that it satisfies 
$S^\circ S^\circ = S^\circ$. 
Let $\varphi \: S^\circ \to B(V)$ be a bounded positive definite 
function with the corresponding representation $\pi_\varphi\colon S^\circ\to B(\cH_\varphi)$ 
on the reproducing kernel Hilbert space $\cH_\varphi$ consisting of $V$-valued functions on $S^\circ$. 
Then the formula 
\begin{equation}
  \label{eq:ext-form}
 (\hat\pi_\varphi(s)f)(t) = f(ts)
\end{equation}
defines a representation $\hat\pi_\varphi \: S \to B(\cH_\varphi)$ by contractions, 
which extends the representation~$\pi_\varphi$. 
\end{prop} 

\begin{prf} Since $\varphi \: S^\circ \to B(V)$ is a bounded positive definite function, 
we obtain the corresponding GNS representation $(\pi_\varphi, \cH_\varphi)$ 
of $S^\circ$, where $\cH_\varphi \subeq V^{S^\circ}$.
The assumption $S^\circ S^\circ = S^\circ$ implies that this 
representation is non-degenerate (Lemma~\ref{L1.1}). 

Since $S$ acts in the obvious way by multipliers on the semigroup 
ideal $S^\circ$, Lemma~\ref{lem:3.13} implies the existence of a bounded 
representation $\hat\pi_\phi$ of $S$ on $\cH_\phi$ which is uniquely 
determined by $\hat\pi_\phi(s)\pi_\phi(t) = \pi_\phi(st)$ for 
$s \in S$, $t \in S^\circ$. 
To see that this representation is also given by \eqref{eq:ext-form},
it suffices to verify this on functions of the form 
$f = \pi_\varphi(t) h$, $h \in \cH_\varphi$: 
\[  (\hat\pi_\varphi(s)\pi_\varphi(t) h)(x) 
=   (\pi_\varphi(st)h)(x) =   h(xst) = (\pi_\varphi(t) h)(xs). \] 
This completes the proof. 
\end{prf}

The following lemma will be a key tool to extend representations from algebra ideals. 

\begin{lem} \mlabel{lem:ideal-lin-ext} Let $\cB$ be a real involutive algebra 
and $\cA \trile \cB$ be an ideal for which $(\cA,p)$ is a seminormed  
involutive algebra. If $\pi \: \cA \to B(\cH)$ is a non-degenerate 
linear $*$-representation bounded on $\ball(\cA,p)$ and 
$S := \{ b \in \cB \: b \ball(\cA,p) \cup \ball(\cA,p)b\subeq \ball(\cA,p)\}$ 
generates~$\cB$ linearly, then there exists a unique linear representation 
$\hat\pi \: \cB \to B(\cH)$ satisfying 
\begin{equation}
  \label{eq:mulrel}
 \hat\pi(b) \pi(a) = \pi(ba) \quad \mbox{ for } \quad a \in \cA, b \in \cB.
\end{equation}
This representation is bounded on $S$. 
\end{lem}

\begin{prf} In view of Lemma~\ref{lem:3.13}, 
the bounded representation  $\pi\res_{\ball(\cA,p)}$ extends to a 
boun\-ded representation  
$\tilde\pi \: S \to B(\cH)$ which is  uniquely determined by the relation 
\[ \tilde\pi(s)\pi(a) = \pi(sa) \quad \mbox{ for }\quad 
s \in S, a \in \ball(\cA,p).\] 

On the other hand, by \cite[Prop.~II.4.14]{Ne00}, there is 
a $*$-representation $\hat\pi \: \cB \to \End(\cH^0)$ on the dense subspace 
$\cH^0 := \Spann \pi(\cA)\cH$ of $\cH$ which is uniquely determined 
by the relation $\hat\pi(b) \pi(a)v = \pi(ba)v$ for 
$a \in \cA, b \in \cB, v \in \cH$. The uniqueness and the linearity of 
$\pi$ now imply that $\hat\pi$ is linear. Since $S$ spans $\cB$ 
and $\tilde\pi(s)\res_{\cH^0} = \hat\pi(s)$ for $s \in S$, it further 
follows that all operators $\hat\pi(b)$, $b \in \cB$, 
are bounded, so that we actually 
obtain a $*$-representation 
$\hat\pi \: \cB \to B(\cH)$ satisfying \eqref{eq:mulrel}. 
\end{prf}

We shall also need the following result. 

\begin{lem}\mlabel{Arv2.4}
Consider the abelian multiplicative semigroup $\cI=((0,1), \cdot)$ 
endowed with the trivial involution. Then the following assertions hold: 
\begin{itemize}
\item[\rm(i)] If $\cB$ is a $C^*$-algebra and $\varphi\colon \cI\to\cB$ is a 
positive definite function, then,  
for every $s\in \cI$, we have $0\le \varphi(s)\in\cB$. 
\item[\rm(ii)] If, moreover, $\cB=\C$ and $\phi$ is bounded, 
then $\varphi$ is also increasing and continuous.
\item[\rm(iii)] Every bounded $*$-representation $(\pi, \cH)$ of $\cI$ is norm-continuous
and extends to a strong\-ly continuous representation $\hat\pi$ 
of the unital semigroup $(0,1]$. 
Moreover, if $\pi$ is non-degenerate and $V \subeq \cH$ is dense, then 
$\pi(r)V$ is dense for every $r \in \cI$, and $\hat\pi(1) = \1$. 
\item[\rm(iv)] 
Let $V$ be a complex Hilbert space and $(A_n)_{n \geq 0}$ in $B(V)$ be a 
sequence for which the series $\varphi(t)=\sum\limits_{n\ge 0}t^nA_n$ is norm convergent in 
$B(V)$ for every $t\in \cI$ 
and the function $\varphi\colon \cI\to B(V)$ is positive definite. 
Then $0\le A_n \le \sup\limits_{0 <  t\le 1}\Vert\varphi(t)\Vert\cdot \1$ 
for every $n\ge 0$. 
\end{itemize}
\end{lem}

\begin{prf} (i) follows from $\varphi(s)=\varphi((s^{1/2})^*s^{1/2})\ge 0$ 
for every $s\in \cI$. 

(ii) From (i) we know that $\varphi(\cI)\subseteq[0,\infty)$. 
Now the conclusion follows by \cite[Prop.~2.4]{Arv87}. 

(iii) Since we may w.l.o.g.\ assume that $\pi$ is non-degenerate and 
$\cI \cong ((0,\infty),+)$, the first assertion 
follows from \cite[Lemma~VI.2.2]{Ne00}. 

Assume that $\pi$ is non-degenerate and that $V \subeq \cH$ is dense. 
If $w \in \cH$ is orthogonal to $\pi(r)V$ for some $r \in \cI$, then 
$\la \pi(r)w, V\ra = \la w, \pi(r)V \ra = \{0\}$, which leads to $\pi(r)w =0$. 
In view of the Spectral Representation $\pi(r) = \int_0^\infty r^\alpha\, dP(\alpha)$ 
(see \cite[Ch.\ X, \S 2, subsect.\ 141]{RSzN72}), 
the fact that $r^\alpha > 0$ for every $\alpha \geq 0$ implies that $\pi(r)$ is injective, 
and therefore $w= 0$. This shows that $\pi(r)V$ is dense for every $r \in \cI$.
We also obtain $\hat\pi(1) = \lim\limits_{r \to 1} \pi(r) = \int_0^\infty \, dP(\alpha) = \1$. 

(iv) For every vector $v\in V$, the function 
$\phi^v \: \cI \to \C, s \mapsto \la\varphi(s)v,v\ra$ 
is positive definite and by (iii) it extends continuously to a positive 
definite function on $[0,1]$. Hence \cite[Lemma~3.8]{Arv87} shows that 
$0\le \la A_nv,v\ra\le\sup\limits_{0\le t\le 1}\la\varphi(t)v,v\ra$ for every $n\ge 0$, 
which implies the assertion. 
\end{prf}

Finally, here is the {\bf proof of Theorem~\ref{thm:Arv2.2}:} 

\begin{prf} We shall prove the dilatability of $\phi$ by verifying 
that, for every $v \in V$, the function $\phi \cdot v$ belongs to  
$\cH_\phi$ (Proposition~\ref{prop:exis-dil}). %{prop:minuni}
Proposition~\ref{P1.2} applies in particular to structured 
$*$-semigroups because the subsemigroup $\cI := \{ \gamma(r) \: 0 < r < 1\}$ of 
$S^\circ$ satisfies 
\[ S^\circ = S \cI = \cI S \quad \mbox{ and }\quad 
S^\circ= S^\circ S^\circ\] 
(cf.\ Remark~\ref{rem:3.16}). 
By construction, the space $\cH :=\cH_\varphi \subeq V^{S^\circ}$ 
consists of functions on $S^\circ$. 
This subsemigroup satisfies 
$\cI \cS^\circ = \cS^\circ = \cS^\circ \cI,$ 
which implies that the representation 
$\pi_\varphi\circ \gamma$ of $(0,1)$ is also non-degenerate. 
From \cite[Lemma~VI.2.2]{Ne00} it now follows that this representation 
of $((0,1),\cdot) \cong ((0,\infty),+)$ is norm-continuous 
and extends to a strongly continuous representation 
of the unital semigroup $(0,1]$. 

For $v \in V$ and $s \in S^\circ$, we have
\begin{equation}
  \label{eq:2}
 (K_s^* v)(t) = \varphi(ts^*)v \quad \mbox{ and } \quad 
\pi_\varphi(s) K_t^* = K_{ts^*}^*. 
\end{equation}
For $v,w \in V$ and $s,t \in S^0$, we thus obtain 
\begin{align*}
 \la \varphi(st^*)v,w \ra 
&= \la K_t^* v, K_s^* w \ra 
= \lim_{r \to 1-} \la \pi_\varphi(r) K_t^* v, K_s^* w \ra  
= \lim_{r \to 1-} \la \varphi(srt^*) v,w \ra.
\end{align*}
Since $S^\circ = S^\circ \cI = S^\circ S^\circ$, this implies that 
\begin{equation}
  \label{eq:e.1}
 \varphi(s) = \lim_{r \to 1-} \varphi(sr) \quad \mbox{ for } \quad 
s \in S^\circ
\end{equation}
in the weak operator topology. 
From Lemma~\ref{Arv2.4} we obtain the relation 
\[ \varphi(\gamma(rr')) \leq \varphi(\gamma(r)) \quad \mbox{ for } \quad 0 < r,r' < 1,\] 
so that the boundedness of $\varphi$ on $\cI$ implies the existence of the strong 
operator limit 
\[ A := \lim_{r \to 1-} \varphi(\gamma(r))  \in B(V). \] 
For $v \in V$ and $0 < r < s < 1$, we further get 
\begin{align*}
\|K^*_{\gamma(r)}v - K^*_{\gamma(s)}v\|^2 
&= \|K^*_{\gamma(r)}v\|^2 + \|K^*_{\gamma(s)}v\|^2 
- 2 \Re \la K^*_{\gamma(r)}v, K^*_{\gamma(s)}v \ra \\
&= \la \varphi(\gamma(r^2)) v,v\ra + \la \varphi(\gamma(s^2)) v,v\ra 
- 2 \Re \la \varphi(\gamma(rs))v, v \ra \to 0 
\end{align*}
for $r,s \to 1-$. This implies that the limit 
\[ \iota(v) := \lim_{r \to 1-} K_{\gamma(r)}^*v = \phi \cdot v\] 
exists for every $v \in V$ in $\cH_\phi$ with 
\[  \|\iota(v)\|^2 
= \lim_{r \to 1-} \|K_{\gamma(r)}^*v\|^2 
= \lim_{r \to 1-} \la K(\gamma(r),\gamma(r))v,v\ra 
= \lim_{r \to 1-} \la \phi(\gamma(r^2))v,v\ra = \la Av,v\ra,\]
which leads to $\|\iota\|^2 \leq \|A\| \leq \limsup\|\phi(\gamma(r))\|$. 
In particular, $\phi \cdot v \in \cH_\phi$ holds for every $v \in V$, 
so that the assertion follows from Proposition~\ref{prop:minuni}. 
\end{prf}

\subsection{Extreme points and homogeneity} 

For later use we also state a rather standard result on the situation when 
the representation $\pi_\vphi$ in Theorem~\ref{thm:Arv2.2} is irreducible 
(compare \cite[Cor. III.1.10]{Ne00} or \cite[Lemma 5.3]{Arv87}). 
We will use the following terminology. 

\begin{defn}\label{homog}
If $S$ is a structured $*$-semigroup and $E$ is a real vector space, 
then a function $f\colon S^\circ\to E$ is said to be  
\emph{homogeneous of degree~$\alpha\in\R$} 
if, for every $x\in S^\circ$ and $0 < r <1$, we have $f(rx)=r^\alpha f(x)$. 
\end{defn}

\begin{prop}\mlabel{Arv5.2}
Let $S$ be a structured $*$-semigroup 
and $\P_1(S^\circ)$ the set of all positive definite functions 
$\varphi\colon S^\circ\to \C$ 
satisfying $\sup\limits_{x\in S^\circ}|\varphi(x)|\le 1$. 
Then $\P_1(S^\circ)$, endowed with the topology of pointwise convergence, 
is a compact convex set. 
For every non-zero extreme point $\varphi \in \P_1(S^\circ)$,  
its corresponding representation $\pi_\phi\colon S^\circ\to C(\cH_\phi)$ 
given by {\rm Theorem~\ref{thm:Arv2.2}} is irreducible and 
there exists some $\alpha\in[0,\infty)$ 
for which $\varphi$ is homogeneous of degree~$\alpha$ and 
\[ \pi_\phi(\gamma(r))= r^\alpha \1 \quad \mbox{ for } \quad 0 < r < 1.\] 
If $\varphi$ is non-constant and $t\mapsto\varphi(tx)$ is continuous on $[0,1)$ 
for all $x\in S^\circ$, then $\alpha>0$.  
\end{prop}

\begin{prf} Let $\bD := \{ z \in \C \: |z| \leq 1\}$. 
Since positive definiteness is preserved under 
pointwise limits, $\P_1(S^\circ)$ is a closed subset of the compact 
set $\bD^{S^\circ}$, hence compact. 

Let $\phi \in \P_1(S^\circ)$ be an extreme point. 
In view of Lemma~\ref{L1.1}, the representation $(\pi_\vphi, \cH_\vphi)$ 
from Theorem~\ref{thm:Arv2.2} is non-degenerate and $v := \iota(1)$ satisfies 
\[ \|v\|^2 = \|\iota^*\iota\| \leq \sup_{0 < r < 1} \|\phi(r)\| \leq 1.\]
Suppose that we have a direct sum decomposition 
$\pi_\vphi = \pi_1 \oplus \pi_2$ with 
$\cH_\vphi = \cH_1 \oplus \cH_2$ and both summands non-zero. 
Then $v = (v_1, v_2) := \iota(1)$ is a cyclic vector with 
$\vphi(s) = \la \pi_\vphi(s)v,v\ra$ for $s \in S$. As $|\vphi(s)| \leq \|v\|^2\leq 1$ 
and $\vphi$ is a non-zero extreme point, we have $\|v\| = 1$. 
The cyclicity of $v$ implies that both components $v_j$ are non-zero. 
Accordingly 
\[ \vphi = \|v_1\|^2\vphi_1  + \|v_2\|^2 \vphi_2\quad \mbox { with }\quad 
\vphi_j(s) :=\frac{\la \pi_j(s)v_j, v_j \ra}{\|v_j\|^2} \]
is a proper convex combination of the two elements 
$\vphi_j \in \P_1(S^\circ)$. This contradiction implies that 
$\pi_\vphi$ is irreducible.  

Now Schur's Lemma implies the existence of a 
function $c\colon(0,1)\to[0,1]$ with $\pi(\gamma(r)) = c(r)\1$ for $0 < r < 1$. 
From \cite[Lemma~VI.2.2]{Ne00} it now follows that $\pi\circ \gamma$ 
is a norm-continuous representation, so that there exists some $\alpha\ge0$ 
with $c(r)=r^\alpha$ for all $0 < r < 1$. We thus obtain 
\[ \vphi(rs) = \iota^* \pi(rs) \iota 
= \iota^* \pi(\gamma(r)) \pi(s) \iota 
= r^\alpha \iota^* \pi(s) \iota 
= r^\alpha \vphi(s) \quad \mbox{ for }  \quad s \in S^\circ, 0 < r < 1.\] 
This proves the first assertion, and the second follows directly.
\end{prf}

The following example shows that Proposition~\ref{Arv5.2} does not generalize 
in any obvious way to operator-valued functions. 
\begin{ex} Let $S := \{ z \in \C \: |z| \leq 1\} =C(\C)$, 
$S^\circ = \{ z \in \C \: |z| < 1\}$ and, for a complex Hilbert space $V$, 
write $\P_1(S^\circ,B(V))$ the set of all positive definite functions 
$\varphi\colon S^\circ\to B(V)$ 
satisfying $\sup\limits_{x\in S^\circ}\Vert\varphi(x)\Vert\le 1$. 
Then $\P_1(S^\circ,B(V))$, endowed with the topology of pointwise weak convergence, 
is a compact convex set. 

The function 
\[ \vphi \: S^\circ \to M_2(\C)\cong B(\C^2), \quad \vphi(z) = \pmat{1 & 0 \\ 0 & z} \] 
is positive definite with 
$\|\vphi(s)\| = 1$ for every $s \in S^\circ$. 
The corresponding space $\cH_\vphi \subeq (\C^2)^{S^\circ}$ is generated by the function 
\[ K_{s,v}(t) = \vphi(ts^*)v = (v_1, ts^*v_2) \quad \mbox{ for } \quad 
v = (v_1, v_2).\] 
It follows that 
\[ \cH_\vphi \cong \C \pmat{1 \\ 0} \oplus \C z \pmat{0 \\ 1}\] 
is $2$-dimensional and multiplicity free. 

If $\vphi = \lambda \vphi_1 + (1-\lambda) \vphi_2$ with 
$0 < \lambda < 1$ and $\vphi_1, \vphi_2 \in \P_1(S^\circ,B(\C^2))$,  
then the functions $\vphi_j$ correspond to subrepresentations, so that 
the only non-trivial way  to decompose $\vphi$ is by 
\[ \lambda \vphi_1(z) =  \pmat{1 & 0 \\ 0 & 0} \quad \mbox{ and } \quad 
(1-\lambda) \vphi_2(z) =  \pmat{0 & 0 \\ 0 & z}.\] 
This leads to $\|\vphi_1(z)\| = \lambda^{-1}$, so that 
$\vphi_1 \not\in \P_1(S^\circ, B(\C^2))$. 
Consequently, $\vphi \in \P_1(S^\circ, B(\C^2))$ is an extreme 
point, the corresponding representation is not irreducible 
and $\vphi$ is not homogeneous.
\end{ex}

\section{Completely positive functions on $\ball(\cA,p)$} 
\mlabel{sec:4}

In this section we first show that a dilation of a positive definite 
function on $\ball(\cA,p)$ is completely positive if and only if $\phi$ is. 
Then we turn to the important special case $\cA = \R$ which provides the 
key tools to obtain a series expansion $\phi = \sum\limits_{n = 0}^\infty \phi_n$ 
of a bounded completely positive function $\phi$ on $\ball(\cA,p)$ into homogeneous 
components $\phi_n$ (Subsection~\ref{subsec:apptoA}). 

\subsection{Complete positivity of dilations} 

The following proposition connects the complete positivity 
of a function $\phi \: \ball(\cA,p) \to B(V)$ to the complete positivity 
of the corresponding GNS representation $(\pi_\phi, \cH_\phi)$. If 
$\cA$ is unital, then $\phi$ is dilatable by Theorem~\ref{thm:Arv2.2}, 
so that it asserts in this case that 
the complete positivity of $\phi$ and $\pi_\phi$ are equivalent. 
It generalizes \cite[Th. 2.17]{Arv87}, which deals with the 
case $\dim V = 1$ and a unital $C^*$-algebra $\cA$. 

\begin{prop}\mlabel{prop:4.1} 
Let $(\cA,p)$ be a seminormed $*$-algebra, $V$  be a complex Hilbert space, 
and $\varphi\colon \ball(\cA,p)\to B(V)$ be a bounded positive definite function. 
Then the following assertions hold: 
\begin{itemize}
\item[\rm(a)] If $\phi$ is completely positive, then the corresponding 
GNS representation $(\pi_\phi, \cH_\phi)$ is completely positive. 
\item[\rm(b)] If $(\pi, \cH,\iota)$ is a minimal dilation of $\phi$, 
then $\varphi$ is completely positive if and only if $\pi$ is completely positive. 
\end{itemize}
\end{prop}

\begin{prf} 
First note that for every $a\in\ball(\cA,p)$ one has $\pi_\phi(a)\in B(\cH_\phi)$ by Remark~\ref{SzN-rem}(a). 

(a) Put $(\pi, \cH) := (\pi_\phi, \cH_\phi)$. 
We have to show that, if $n\ge 1$, $A := (a_{jk})\in M_n(\cA)_+$ 
and $a_{jk}\in \ball(\cA,p)$ for $1\le j,k\le n$,  
then $(\pi(a_{jk}))\in M_n(B(\cH))_+$. 
Since the functions $K_s^*v$, $s \in \ball(\cA,p)$, $v \in V$, span a dense 
subspace of $\cH$, 
it suffices to prove that 
$\sum\limits_{j,k=1}^n\la\pi(a_{jk})h_k,h_j\ra\ge0$ 
for arbitrary $h_1,\dots,h_n\in \Spann \{  K_s^* v \: s \in \ball(\cA,p), v \in V\}$. 

There exists an $m \in \N$, 
and $b_{jp}\in \ball(\cA,p)$ and $v_{jp}\in V$ for $1\le j\le n$ and $1\le p\le m$, 
such that $h_j=\sum\limits_{p=1}^m K_{b_{jp}}^*v_{jp}$. 
Then 
\begin{align*}
&\sum_{j,k=1}^n\la\pi(a_{jk})h_k,h_j\ra
=\sum_{j,k=1}^n\sum_{p,q=1}^m\la\pi(a_{jk})K_{b_{kp}}^*v_{kp},K_{b_{jq}}^*v_{jq}\ra \\
&=\sum_{j,k=1}^n\sum_{p,q=1}^m\la K_{b_{jq}} \pi(a_{jk}) K_{b_{kp}}^*v_{kp},v_{jq}\ra 
=\sum_{j,k=1}^n\sum_{p,q=1}^m\la\varphi(b_{jq}a_{jk}b_{kp}^*)v_{kp},v_{jq}\ra,
\end{align*}
and since $\varphi$ is completely positive, it suffices to check that 
the matrix 
$$C:=(b_{jq}a_{jk}b_{kp}^*)\in M_{mn}(\cA)$$ 
is non-negative. 
If we denote by $B_p\in M_n(\cA)$ the diagonal matrix 
$\diag(b_{1p},\dots,b_{np})$, then 
$C = D (A \otimes \1_m) D^*\ge 0$, 
where 
$$D=\begin{pmatrix}
B_1 & 0 & \dots & 0  \\
B_2 & 0 & \dots & 0 \\
\vdots & & \vdots \\
B_m & 0 & \dots & 0
\end{pmatrix}\in M_m(M_n(\cA))\cong M_{mn}(\cA) $$
and this proves that $\pi$ is completely positive. 

(b) It is well known (and easy to see) that the linear mapping 
$\Psi\colon B(\cH)\to B(V)$, $T\mapsto \iota^*T\iota$ is completely positive. 
Therefore $\phi = \Psi \circ \pi$ is completely positive if $\pi$ is completely positive. 
For the converse we can use (a) because  $(\pi, \cH)$ is equivalent to the 
GNS representation $(\pi_\phi, \cH_\phi)$ by Proposition~\ref{prop:minuni}.
\end{prf}

\begin{cor} \mlabel{cor:4.2} A bounded $*$-representation $(\pi, \cH)$ of 
$\ball(\cA,p)$ is completely positive if and only if, for 
every $v$ in some dense subset of~$\cH$, the function 
$\pi^v(a) := \la \pi(a)v,v\ra$ is completely positive. 
\end{cor}

\begin{prf} If $\pi$ is completely positive, then so is every function 
$\pi^v$ since the linear functional $A \mapsto \la Av,v\ra$ on $B(\cH)$ is 
positive definite, hence completely positive (Remark~\ref{rem:1.5}(c)). 

Suppose, conversely, that all functions $\pi^v$ are completely positive 
for $v$ in some dense subset of $\cH$, hence actually for every $v\in\cH$. 
Write $\pi$ as the direct sum of a zero representation and a 
direct sum of cyclic representations $(\pi_j, \cH_j,v_j)_{j \in J}$. 
Since $v_j\in\cH_j$ is a cyclic vector for the representation $\pi_j$, 
it follows that $\pi_j$ is a minimal dilation of the positive definite function $\phi_j:=\pi^{v_j}$, 
and then $\pi_j$ is unitarily equivalent to the GNS representation $\pi_{\phi_j}$ by Proposition~\ref{prop:minuni}(b). 
On the other hand, since the function $\phi_j$ is completely positive by assumption, 
the representation $\pi_{\phi_j}$ is completely positive by Proposition~\ref{prop:4.1}(a). 
Hence also $\pi_j$ is completely positive. 
 This implies that $\pi$ is completely 
positive, as an orthogonal direct sum of completely positive maps. 
\end{prf}

\begin{cor}\label{Arv2.17_cor}
Let $(\cA,p)$ be a seminormed $*$-algebra and $V$ be a complex Hilbert space.  
If $\varphi\colon \ball(\cA,p)\to B(V)$ is a dilatable completely positive function, 
then $\varphi-\varphi(0)$ is also completely positive and dilatable. 
\end{cor}

\begin{prf} Let $(\pi, \cH,\iota)$ be a dilation of $\phi$. 
With the same argument as in the proof of Lemma~\ref{Arv2.2-cor2}, 
we see that $P := \pi_\phi(0)$ is a projection onto a trivial subrepresentation. 
With $\cH_0 := P(\cH)$ and $\cH_1 := \cH_0^\bot$ and the corresponding decompositions 
$\iota = \iota_0  + \iota_1$ and $\pi(s) = (\pi_0(s), \pi_1(s)) = (\1,\pi_1(s))$, we then obtain 
\[ \phi(s) = \iota^* \pi(s) \iota = \iota_0^* \iota_0 + \iota_1^* \pi_1(s) \iota_1 
\quad \mbox{ and } \quad 
\phi(0) = \iota^* \pi(0) \iota = \iota_0^* \iota_0.\] 
This proves the assertion. 
\end{prf}

\begin{lem} \mlabel{lem:4.4} 
Let $(\pi, \cH)$ be a non-degenerate bounded representation of $\ball(\cA,p)$, 
with its canonical extension 
$\pi^1 \: \ball(\cA^1, p^1) \to B(\cH)$ 
{\rm(Example~\ref{ex:3.14})}. 
Then  $\pi$ is completely positive if and only if $\pi^1$ is completely positive. 
\end{lem}

\begin{prf} We must show that the complete positivity of $\pi$ implies the complete 
positivity of $\pi^1$. Decomposing $\pi$ as a direct sum of cyclic representations, 
it suffices to prove the lemma for cyclic representation. So assume that 
$0 \not= v \in \cH$ is cyclic and put 
$\phi(a) := \la \pi(a)v,v\ra$. For $s_1, \ldots, s_n \in \ball(\cA,p)$ and 
$c_1, \ldots, c_n \in \C$, the vectors of the form 
$w := \sum_{j = 1}^n c_j \pi(s_j)v$ form a dense subspace of $\cH$. 
In view of 
Corollary~\ref{cor:4.2}, it therefore suffices to show that the 
functions 
\[ \psi \: \ball(\cA^1, p^1) \to \C, \quad 
\psi(a) := \la \pi^1(a)w, w \ra 
= \sum_{i,j=1}^n c_i \oline{c_j} \la \pi(s_j^* a s_i) v, v \ra 
= \sum_{i,j=1}^n c_i \oline{c_j} \phi(s_j^* a s_i) \] 
are completely positive. 

So let $A = (a_{k\ell}) \in M_N(\cA^1)$ be positive. Then 
the matrix $B := (s_j^* a_{k\ell} s_i) \in M_{nN}(\cA)$ is positive by 
Example~\ref{ex:2.17}(a). Now the complete positivity of $\phi$ implies that the 
matrix $(\phi(s_j^* a_{k\ell} s_j)) \in M_{nN}(\C)$ is positive, so that 
Example~\ref{ex:2.17}(b) shows that the matrix 
\[ \psi_N(A) = \Big(\sum_{i,j=1}^n c_i \oline{c_j} \phi(s_j^* a_{k\ell} s_i)\Big)_{1 \leq k,\ell \leq N} 
\in M_N(\C) \] 
is positive. We conclude that $\psi$ is completely positive. 
\end{prf}

The following proposition is a version of Proposition~\ref{Arv5.2} for completely positive functions 
on $\ball(\cA,p)$. 

\begin{prop}\mlabel{Arv5.2-compl}
If $(\cA,p)$ is a seminormed involutive algebra,  
then the set $\CP_1(\ball(\cA,p))$ of all 
completely positive functions $\varphi\colon \ball(\cA,p)\to \C$ with 
\[ \|\phi\|_\infty := \sup\limits_{x\in \ball(\cA,p)}|\varphi(x)|\le 1,\] 
endowed with the topology of pointwise weak convergence, is a compact convex set. 
For every extreme point $\varphi \in \CP_1(\ball(\cA,p))$ 
either $\cH_\phi = \{0\}$ 
or the representation $\pi_\phi\colon S^\circ\to C(\cH_\phi)$  
%given by {\rm Theorem~\ref{thm:Arv2.2}} 
is irreducible and there exists some integer $\alpha \in \N_0$ such that 
\begin{equation}
  \label{eq:homog}
\phi(tab) = t^\alpha \phi(ab) \quad \mbox{ for } \quad 0 \leq |t| < 1, 
a,b \in \ball(\cA,p).
\end{equation}
\end{prop}

\begin{prf} Note that complete positivity is preserved under 
pointwise convergence, so that the same argument as in Proposition~\ref{Arv5.2} 
implies the compactness of $\CP_1(\ball(\cA,p))$. 
We likewise see that, for every extreme point 
$\vphi$ of $\CP_1(\ball(\cA,p))$, for which the corresponding Hilbert space 
$\cH_\phi$ is non-zero, the corresponding representation 
$\pi_\vphi$ is irreducible. Here we use that the complete positivity of 
$\phi$ implies the complete positivity  of $\pi_\phi$ (Proposition~\ref{prop:4.1}(a)), 
and therefore the complete positivity of the functions 
$\vphi_j(s) :=\frac{\la \pi_j(s)v_j, v_j \ra}{\|v_j\|^2}$. 
Using the canonical extension of the representation to 
$\ball(\cA^1,p^1)$ (Example~\ref{ex:3.14}) and Schur's Lemma, 
this implies the existence of some $\alpha \geq 0$ with 
$\tilde\pi_\phi(t) = t^\alpha \1$ for $|t| <1$. 
For $v \in V$ and $a \in \ball(\cA,p)$, we now have
\[ t^\alpha K_{a,v}(b) = (\tilde\pi_\phi(t) K_{a,v})(b) 
= K_{a,v}(tb) = \phi(tba^*)v,\] 
so that $\phi$ satisfies \eqref{eq:homog}. 
It remains to show that $\alpha \in \N_0$. 

If $\la \phi(a^*a)v,v\ra =0$ for every $a \in \ball(\cA,p)$ and $v \in V$, 
then the relation $\|K_{a,v}\|^2 = \la \phi(a^*a)v,v\ra$ in $\cH_\phi$ implies 
that $\cH_\phi = \{0\}$. If this is not the case, then there exists an 
$a \in \ball(\cA,p)$ and $v \in V$ such that the completely positive function 
$\psi(t) := \la \phi(a^*ta)v,v\ra = t^\alpha \la \phi(a^*a)v,v\ra$ on $[-1,1]$ does 
not vanish. Its complete positivity implies that, for every positive definite matrix  
$A = (a_{ij}) \in M_n(\R)$ with $a_{ij} > 0$, the matrix 
$(a_{ij}^\alpha) \in M_n(\R)$ is positive semidefinite. 
In view of \cite[pp.~270-271]{Ho69}, this implies that $\alpha \in \N_0$ 
(see also \cite[Thm.~1]{CR78}). 
\end{prf}

\begin{prop} \mlabel{prop:4.3} Let $(\cA,p)$ be a seminormed involutive algebra and 
$\phi \: \cA \to B(V)$  be dilatable 
completely positive and homogeneous of degree $n \in \N_0$ and 
bounded on $\ball(\cA,p)$. Then $\phi$ is a continuous polynomial which is 
homogeneous of degree $n$. 
\end{prop}

\begin{prf} {\bf Case 1:} First we assume that $\cA$ is unital. We 
use the method of proof of \cite[Lemma~4.3]{Arv87} to see that, for $h \in \cA$, the function 
\[ (\Delta_h \phi)(a) := \phi(a+ h) - \phi(a) \] 
is also completely positive. Iterating this argument, one shows that the function 
\[ (\Delta^n \phi)(a; h_1, \ldots, h_n) := 
(\Delta_{h_1}\cdots \Delta_{h_n})(a) \] 
is completely positive on the seminormed involutive algebra 
$\cA^{n+1}$ with respect to the seminorm 
\[ q_{n+1}(a_0, \ldots, a_n) := \max \{ p(a_0), \cdots, p(a_n)\}.\] 
The corresponding argument in the proof of \cite[Lemma~5.8]{Arv87} works for general 
unital involutive algebras, 
but one should replace the reference to \cite[Thm.~2.17]{Arv87} in its proof by a reference to 
\cite[Lemma~4.3]{Arv87}. 

Theorem~\ref{thm:Arv2.2} shows that, if $\phi(r) = 0$ for $0 < r < 1$, 
then $\iota = 0$, and thus $\phi =0$. 
Now the same argument as in \cite[p.~342]{Arv87} implies that 
$\phi$ is a homogeneous polynomial of degree $n$. In the whole argument, continuity 
is never needed, only that the functions $\Delta^k \phi$ are bounded on balls with respect to 
the corresponding seminorm $q_{k+1}$. 

{\bf  Case 2:} Now we assume that $\cA$ is not unital. 
Let $(\pi, \cH, \iota)$ be a dilation of $\phi$. 
Since we may assume that $\cH = \cH_\phi$ (Proposition~\ref{prop:minuni}) 
and $\phi$ is homogeneous of degree 
$n$, it follows that $\pi(ta) = t^n \pi(a)$ for $a \in \cA$ and $t \in \R$. 
In view of Lemma~\ref{lem:4.4}, the bounded representation 
$\pi$ of $\ball(\cA,p)$ extends to a bounded representation 
$\pi^1$ of $\ball(\cA^1, p^1)$. Further Proposition~\ref{prop:4.1} implies that 
$\pi^1$ is completely positive. Therefore 
the first part of the proof applies to the completely positive function 
$\phi^1(a) := \iota^* \pi^1(a) \iota$ which extends $\phi$ to $\ball(\cA^1, p^1)$. 
\end{prf}

Since representations are dilatable for trivial reasons (Example~\ref{ex:3.8}), 
the preceding proposition implies: 

\begin{cor} If $\pi \:  (\cA,\cdot) \to B(\cH)$ 
is a representation bounded on $\ball(\cA,p)$ and 
homogeneous of degree $n$, then $\pi$ is a continuous homogeneous polynomial 
of degree $n$. 
\end{cor}

\begin{rem} Let $\phi \:  \cA \to B(V)$  be completely positive, 
homogeneous of degree $n \in \N_0$ and 
bounded on $\ball(\cA,p)$ but not necessarily dilatable. Then all functions in 
$\cH_\phi \subeq \C^{\ball(\cA,p)}$ are homogeneous of degree $n \in \N_0$, 
which implies that the representation $\pi_\phi$ is homogeneous 
of degree $n$. Now the preceding corollary implies that 
$\pi_\phi$ is a continuous homogeneous polynomial of degree~$n$.
\end{rem}

\begin{rem}
It follows by \cite[Th.~4]{AC86} 
that if $\cA$ is a $C^*$-algebra and $S=(\cA,\cdot)$ its multiplicative group, 
then the completely positive maps $\varphi\colon S\to B(V)$ 
which are homogeneous of degree~1 
in the sense of the above Definition~\ref{homog} 
are precisely the (completely positive) $\R$-linear maps.
\end{rem}

\subsection{The case $\cA = \R$} 
\mlabel{subsec:4.2} 

In this section we take a closer look at the special case $\cA = \R$. We write 
$S := \ball(\R) = ((-1,1), \cdot)$. 
We know from Proposition~\ref{Arv5.2-compl} that the 
set $\CP_1(S)$ of completely positive functions $\phi$ on $S$ 
with $\|\phi\|_\infty \leq 1$ is a compact convex set with respect to the 
topology of pointwise convergence and that its non-zero extreme points 
are the characters 
\[ \chi_n(r) = r^n, \quad n \in \N_0.\] 

Then 
\[ \bS := \Big\{ (a_n)_{n\in \N_0} \in \ell^1(\N_0) \: a_j \geq 0, \sum_n a_n \leq 1\Big\}\] 
is a weak-$*$-closed bounded subset of $\ell^1(\N_0) \cong c_0(\N_0)'$, 
hence compact in the weak-$*$-topology. This implies that on $\bS$ the 
weak-$*$-topology coincides with the product topology inherited from $\R^{\N_0}$. We have a 
well-defined map 
\[ \Phi \: \bS \to \CP_1((-1,1)), \quad 
\Phi((a_n))(t) := \sum_{n = 0}^\infty a_n t^n, \] 
and since the sequences $(t^n)_{n \in \N_0}$, $|t| < 1$, are in $c_0(\N_0)$, 
$\Phi$ is continuous with respect to the product topology on $\bS$ and 
the topology of pointwise convergence on $\CP_1((-1,1))$. 
Since $\Phi$ is affine, its image is a compact convex subset of
$\CP_1((-1,1))$ which contains the extreme points 
$\chi_n = \Phi(e_n)$. Hence $\Phi$ is surjective by the Krein--Milman 
Theorem. As $\Phi$ is obviously injective, 
the compactness of $\bS$ implies that $\Phi$ is a homeomorphism. Hence 
every function $\phi \in \CP_1((-1,1))$ is analytic 
and its Taylor series in $0$ converges uniformly on the closed interval $[-1,1]$. 
In particular, $\phi$ extends to a continuous function on $[-1,1]$. 

\begin{rem} The completely positive characters of $(-1,1)$ are 
the monomials $\chi_n(t) = t^n$, $n \in \N_0$. 
\end{rem}

\begin{prop} \mlabel{prop:4.1b} {\rm(Completely positive functions on $(-1,1)$).} 
  \begin{itemize}
  \item[\rm(a)] Every non-degenerate 
bounded completely positive representation %\break 
$(\pi, \cH)$ of $(-1,1)$ 
is a direct sum of subrepresentations $(\pi_n, \cH_n)_{n \in \N_0}$ with 
$\pi_n(t) = t^n \1$ for $|t| < 1$. 
  \item[\rm(b)] Every bounded completely positive function $\phi \: (-1,1) \to B(V)$ 
is of the form 
\[ \phi(t) = \sum_{n = 0}^\infty t^n A_n,\] where $A_n \geq 0$ and 
the sum converges on each interval 
$[-r,r]$, $r < 1$, in the norm topology. 
Further, $\phi(\pm 1) := \sum_{n = 0}^\infty (\pm 1)^n A_n$ defines a weakly continuous 
extension to $[-1,1]$. 
\end{itemize}
\end{prop}

\begin{prf} (a) Since every non-degenerate representation is a direct sum of cyclic ones, 
this follows from the structure of the bounded completely positive functions 
on $(-1,1)$: 
\[ \phi(t) = \sum_{n = 0}^\infty a_n t^n\quad \mbox{ with } \quad a_n \geq 0 \quad \mbox{ and } \quad 
 \sum_n a_n < \infty.\] 
It implies that $\cH_\phi$ contains $\chi_n$ whenever $a_n > 0$ because 
$\phi - a_n \chi_n$ is positive definite 
(\cite[Thm.~I.2.8]{Ne00}). Let $\phi_n := a_n \chi_n$. Then the norm of 
$\phi_n$ in $\cH_{\phi_n}$ is given by $\|\phi_n\|^2 = a_n$, so that 
the series $\sum_n \phi_n$ converges in $\oplus_{n = 0}^\infty \cH_{\phi_n}$. 
This shows that $\cH_\phi \cong \oplus_{n = 0}^\infty \cH_{\phi_n}$. 

(b) Let $(\pi, \cH, \iota)$ be a dilation of $\phi$ 
(Theorem~\ref{thm:Arv2.2}). Then Proposition~\ref{prop:4.1} implies that the 
non-degenerate representation $\pi$ is completely positive. 
Now (a) implies that 
$(\pi, \cH) = \oplus_{n = 0}^\infty (\pi_n, \cH_n)$ and we obtain with 
$\iota(v) = \sum_n \iota_n(v)$, $\iota_n(v) \in \cH_n$, the relations 
\[ \pi(t) = \sum_n t^n \id_{\cH_n} \quad \mbox{ and }  \quad 
\phi(t) = \sum_n t^n A_n \quad \mbox{ with } \quad A_n := \iota_n^* \iota_n\in B(V),\] 
where the series converge in the weak operator topology and uniformly on $[-r,r]$ in the norm 
topology. We further have the weak operator convergent series 
\[ \iota^*\iota = \sum_n \iota_n^* \iota_n = \sum_n A_n \] 
which leads to weak operator continuous extension of $\phi$ by 
$\phi(\pm 1) := \sum_{n = 0}^\infty (\pm 1)^n A_n$. 
\end{prf}

\subsection{Series expansion of bounded completely positive functions} 
\mlabel{subsec:apptoA}

For a bounded function $\phi \: \ball(\cA,p) \to E$ with values in a Banach space $E$, 
we write 
\[ \|\phi\|_\infty := \sup \{ \|\phi(a)\| \: a \in \ball(\cA,p)\}.\] 

\begin{thm} \mlabel{thm:4.7} 
Let $\phi \: \ball(\cA,p) \to B(V)$ be a bounded dilatable completely positive function. 
Then $\phi$ is analytic in the $p$-topology and its Taylor series 
$\phi = \sum_n \phi_n$ at $0$ 
converges for every $r < 1$ uniformly on the ball 
$\{ a \in \cA \: p(a) \leq r\}$. 
The homogeneous polynomials $\phi_n$ are also completely positive with 
\[ \|\phi_n\|_\infty \leq \|\phi\|_\infty \quad \mbox{ for every } \quad n \in \N.\] 
\end{thm}

\begin{prf} Let $(\pi, \cH,\iota)$ be a 
be a minimal dilation of $\phi$ 
(Theorem~\ref{thm:Arv2.2}). 
Then the representation $\pi$ is non-degenerate and is completely positive by 
Proposition~\ref{prop:4.1}, and by  
Lemma~\ref{lem:4.4} it extends to a completely positive representation 
$\pi^1$ of $\ball(\cA^1,p^1)$. 
Applying Proposition~\ref{prop:4.1b} to 
$\pi^1\res_{(-1,1)\1}$, we see that 
$(\pi,\cH) \cong \bigoplus\limits_{n \in \N_0} (\pi_n, \cH_n)$ with the $\ball(\cA,p)$-invariant subspaces 
\[ \cH_n = \{ v \in \cH \: \pi^1(t\1)v = t^n v \quad \mbox{ for } \quad |t| < 1\}.\] 
Writing $\iota(v) = \sum\limits_n \iota_n(v)$ with $\iota_n(v) \in \cH_n$, this leads to the weakly 
convergent series  
\[ \phi(a) 
= \iota^* \pi(a) \iota 
= \sum_n \phi_n(a) \quad \mbox{ with } \quad 
\phi_n(a) := \iota_n^* \pi_n(a) \iota_n,\] 
where $\pi_n$ and $\phi_n$ are homogeneous of degree $n$ in the sense that 
\[ \phi_n(ta) = t^n \phi_n(a) \quad \mbox{ and } \quad 
\pi_n(ta) = t^n \pi_n(a) \quad \mbox{ for } \quad a \in \ball(\cA,p), |t| < 1.\] 
Since $\|\iota_n\| \leq \|\iota\|$ and $\|\pi_n(a)\| \leq \|\pi(a) \| 
\leq 1$, the series 
$\phi(a)  = \sum\limits_{n = 0}^\infty \phi_n(a)$ 
converges uniformly on all subsets 
$B_r := \{ a \in \cA \: p(a) < r\}$ for $r < 1$ with respect to the 
norm topology on $B(V)$, and pointwise on all of $\ball(\cA,p)$. 
To see that $\phi$ is analytic, it now suffices to observe that each $\phi_n$ is polynomial, 
which follows from Lemmas~\ref{lem:1.5} and Proposition~\ref{prop:4.3}. 
\end{prf}

The main point of the preceding theorem is the fact that 
the homogeneous functions $\phi_n$ are polynomials. This would follow 
much more directly if $\phi$ is weakly smooth, 
because then the Taylor series of the functions 
$a \mapsto \la \phi(a)v,w\ra$, $v,w \in V$, makes sense, 
and the terms of order $n$ are polynomial in~$a$. 

\begin{rem} (On the non-unital case) 
Every completely positive representation \break 
$\pi \: \ball(\cA,p) \to C(\cH)$ is in particular a dilatable completely 
positive function, so that Theorem~\ref{thm:4.7} applies. 

This applies in particular to the GNS representation 
$(\pi,\cH)= (\pi_\phi, \cH_{\phi})$ associated 
to a bounded completely positive function $\phi \: \ball(\cA,p) \to B(V)$ 
(Proposition~\ref{prop:4.1} and Remark~\ref{SzN-rem}). 
From the decomposition $\pi = \sum_{n = 0}^\infty \pi_n$, 
we derive that the corresponding positive definite kernel 
$K(a,b) := \phi(ab^*)$ decomposes as 
$K = \sum_{n = 0}^\infty K^n,$ and these kernels satisfy 
\[ K^n(ab,c) = K^n(a,cb^*) \quad \mbox{ for } \quad a,b,c \in \ball(\cA,p).\] 
If $\cA$ is not unital, then it is not clear if there exists functions 
$\phi_n \: \ball(\cA,p)\to B(V)$, $n \in \N_0$, with $\phi_n(ab^*) = K^n(a,b)$ for $a,b \in 
\ball(\cA,p)$. Without the assumption of dilatability it is therefore not clear if 
$\phi$ admits a corresponding series expansion. 
If $\ball(\cA,p)$ contains an approximate identity and $\phi$ is continuous, 
then we can use \cite[Prop.~IV.1.29]{Ne00} to obtain functions 
$\phi_n$ on the dense subsemigroup $\ball(\cA,p)\ball(\cA,p)$ of $\ball(\cA,p)$ which 
satisfy $\phi_n(ab^*) = K_n(a,b)$ for $a,b \in \ball(\cA,p)$. 
\end{rem}

\section{Analytic positive definite functions} 
\mlabel{sec:5}

In this section we turn to positive definite functions 
$\phi \: \ball(\cA,p) \to B(V)$ which are analytic with respect to some 
locally convex topology on $\cA$ for which $p$ is a continuous seminorm 
(cf. Definition~\ref{def:0.2}). After 
some general observations concerning series expansions on structured semigroup 
in Subsection~\ref{subsec:6.1}, we derive in Subsection~\ref{subsec:6.2} a series 
expansion of dilatable bounded analytic positive definite functions 
$\phi \: \ball(\cA,p) \to B(V)$. Before we can use these result to show that 
$\phi$ is completely positive, we take in Subsection~\ref{subsec:6.3} a detailed look 
at the case where $\phi$ is a homogeneous polynomial. The main result to deal 
with this case is Proposition~\ref{prop:sn-univ} which asserts that 
the $C^*$-algebra $S^n(C^*(\cA,p))$ is the universal enveloping algebra of 
$(S^n(\cA),p_n)$. This in turn is used to show that every dilatable positive definite 
homogeneous polynomial $\phi \: \cA \to B(V)$ which is bounded on $\ball(\cA,p)$ is completely 
positive. 

\subsection{Series expansions on structured $*$-semigroups}
\mlabel{subsec:6.1}

\begin{prop}\mlabel{series}
Let $S$ be a structured $*$-semigroup, $T \subeq S$ be a subsemigroup 
with $rT \subeq T$ for $0 < r < 1$ and 
$V$ be a complex Hilbert space. 
For every integer $n\ge 0$, let 
the function $\varphi_n\colon T\to B(V)$ be homogeneous of degree~$n$ 
such that, for every $a\in T$, the series 
\[ \varphi(a):=\sum\limits_{n\ge 0}\varphi_n(a)\] 
is convergent in the weak operator topology 
and $M:=\sup\limits_{a\in T}\Vert\varphi(aa^*)\Vert<\infty$. 
Then the function $\varphi$ is positive definite 
if and only if all the functions $(\varphi_n)_{n \geq 0}$ are positive definite on $T$. 
If this is the case, then 
$\sup\limits_{a\in T}\Vert\varphi_n(aa^*)\Vert\le M$ for all $n\ge 0$.
\end{prop}

\begin{prf}
Let us fix arbitrary $m\ge 1$, $a_1,\dots,a_m\in T$, and $v_1,\dots,v_m\in V$,  
and set $\alpha_n:=\sum\limits_{j,k=1}^m\la\varphi_n(a_j^*a_k)v_k,v_j\ra\in\C$. 
If $0\le t\le 1$, the series 
$\psi(t)=\sum\limits_{n\ge 0}\alpha_nt^n
=\sum\limits_{j,k=1}^m\la\varphi(ta_j^*a_k)v_k,v_j\ra$ 
is convergent. 
Moreover, if $0\le t_1,\dots,t_q\le 1$, and $\mu_1,\dots,\mu_q\in\C$, then 
$$\sum_{i,\ell=1}^q\psi(t_it_\ell)\mu_i\oline\mu_\ell
=\sum_{i,\ell=1}^q\sum\limits_{j,k=1}^m\la\varphi((t_\ell a_j)^*(t_ia_k))(\mu_iv_k),\mu_\ell v_j\ra\ge0$$
since $\varphi$ is positive definite. 
Therefore $\psi\colon[0,1]\to\C$ is positive definite on the $*$-semigroup $[0,1]$, 
and then  \cite[Lemma 3.8]{Arv87} shows that, for every $n\ge 0$, we have 
\[ 0\le\alpha_n\le \sup\limits_{0\le t\le 1}\psi(t),\]  
that is 
\begin{equation}\label{Arv3.9}
0\le 
\sum\limits_{j,k=1}^m\la\varphi_n(a_j^*a_k)v_k,v_j\ra
\le\sup\limits_{0\le t\le 1}\sum\limits_{j,k=1}^m\la\varphi(ta_j^*a_k)v_k,v_j\ra.
\end{equation}
Recalling the definition of $\alpha_n$, it follows that, for each $n\ge 0$, 
the homogeneous function $\varphi_n\colon T \to B(V)$ is positive definite. 
Using %Lemma~\ref{Arv_under3.9} {\bf replace ref!} and 
the above equation~\eqref{Arv3.9} (with $m=1$),  
we obtain 
\[ \sup\limits_{a\in T}\Vert\varphi_n(a^*a)\Vert
\le \sup\limits_{a\in T}\Vert\varphi(a^*a)\Vert = M.\] 
This concludes the proof.
\end{prf}

\begin{cor}\mlabel{series-2}
Let $S$ be a structured $*$-semigroup and $V$ be any complex Hilbert space. 
For every integer $n\ge 0$, let 
the function $\varphi_n\colon S^\circ\to B(V)$ be homogeneous of degree~$n$ 
such that, for every $a\in S^\circ$, the series 
\[ \varphi(a):=\sum\limits_{n\ge 0}\varphi_n(a)\]
is convergent in the weak operator topology 
and $M:=\sup\limits_{a\in S^\circ}\Vert\varphi(a)\Vert<\infty$. 
Then the function $\varphi$ is positive definite 
if and only if all the functions $(\varphi_n)_{n \geq 0}$ are positive definite. 
If this is the case, then 
$\sup\limits_{a\in S^\circ}\Vert\varphi_n(a)\Vert\le M$ for all $n\ge 0$.
\end{cor}

\begin{prf} 
We apply Proposition~\ref{series} to $T = S^\circ$ and 
recall that $S^\circ = S^\circ S^\circ$, which implies that, for 
$s = ab^*$, we obtain 
\[ \|\phi_n(s)\|^2 = \|\phi_n(ab^*)\|^2 \leq \|\phi_n(aa^*)\|\|\phi_n(bb^*)\| \leq M^2.
%\qedhere
\] 
This completes the proof. 
\end{prf}

\subsection{Analytic positive definite functions} 
\mlabel{subsec:6.2}

\begin{defn} If $M$ is an analytic manifold and $V$ a complex Hilbert space, 
then we call a function $\phi \: M \to B(V)$ {\it weakly analytic} if there exists a 
dense subspace $E \subeq V$, such that, the function 
$M \to \C, m \mapsto \la \phi(m)v,w\ra$ is analytic for $v,w \in E$. 
\end{defn}

\begin{prop} \mlabel{prop:5.7} 
Every bounded weakly analytic positive definite function $\phi \: (-1,1) \to B(V)$ is completely positive. 
\end{prop}

\begin{prf} (a) First we consider the case $V = \C$. 
Then $\phi\res_{[0,1)}$ is increasing by Lemma~\ref{Arv2.4}(ii), 
and by \cite[Cor.\ VI.2.11]{Ne00} there exists a 
uniquely determined positive Radon measure $\mu$ on $[0,\infty)$ with 
\[ \phi(t) = \int_0^\infty t^\alpha\, d\mu(\alpha) \quad \mbox{ for } \quad 0 < t < 1.\]  On the other hand, the analyticity of $\phi$ in $0$ implies the existence of an $\eps > 0$ and 
$a_n \in \C$ with 
\[ \phi(t) = \sum_{n = 0}^\infty a_n t^n \quad \mbox{ for } \quad |t| \leq \eps.\] 
Then the function $\psi(t) := \phi(t\eps) = \sum_{n = 0}^\infty a_n \eps^n t^n$ on $[0,1]$ is 
is positive definite, so that Lemma~\ref{Arv2.4}(iv) implies that 
$0 \leq a_n \leq \eps^{-n} \psi(1) 
= \eps^{-n}\phi(\eps)$. 
From the uniqueness of the measure $\mu$ and the identity 
\[ \int_0^\infty t^\alpha\, \eps^\alpha d\mu(\alpha) = \phi(t\eps) = \sum_{n = 0}^\infty a_n \eps^n t^n \quad \mbox{ for } \quad 0 < t \leq 1, \] 
 it follows that $\mu = \sum_{n = 0}^\infty a_n \delta_n$. This in turn implies that 
\[ \phi(t) = \int_0^\infty t^\alpha\, d\mu(\alpha) = \sum_{n = 0}^\infty a_n t^n \] 
holds for $0 \leq t <1$ and by analyticity we have this relation for $|t| < 1$. 
It follows in particular that $\phi$ is completely positive. 

(b) Now let $\phi \: (-1,1) \to B(V)$ be weakly analytic and positive definite. 
Let $E \subeq V$ be a dense subspace such that $s \mapsto \la \phi(s)v,w \ra$ is analytic for 
$v,w \in E$. 
To see that $\phi$ is completely positive, it suffices to show that the corresponding 
GNS representation $\pi_\phi$ is completely positive (Proposition~\ref{prop:4.1}). 
This representation is weakly analytic, because the functions 
\[s\mapsto 
\la \pi_\phi(s) \pi_\phi(t) \iota(v), \pi_\phi(u)\iota(w) \ra = \la \phi(ust)v,w \ra,  
\qquad v,w \in E, \] 
are analytic for every $t, u \in (-1,1)$. 
Therefore we may assume that $\phi = \pi$ is a non-degenerate 
representation. 
Then $\pi$ is a direct sum of cyclic representations whose generating vectors 
have analytic matrix coefficients. 
Combining (a) with Proposition~\ref{prop:4.1}, it follows that 
all these subrepresentations are completely positive, and therefore $\pi$ is completely positive. 
\end{prf}

Let $\vphi \: S  := \ball(\cA,p) \to B(V)$ be a bounded positive 
definite real analytic function for some locally convex topology for which 
$p$ is continuous (Definition~\ref{def:0.2}(c)) and 
$\pi_\vphi \: S \to B(\cH_\phi)$ be the corresponding contraction representation 
of $S$ (Theorem~\ref{thm:Arv2.2}). 
Since $\phi$ is analytic, there exists an open circular $0$-neighborhood 
$U \subeq S$ such that 
\[ \phi(a) = \sum_{n \geq 0} \phi_n(a) \quad \mbox{ with } \quad 
\phi_n(a,\ldots, a) = \frac{1}{n!} (\partial_a^n \phi)(0)\] 
converges for every $a \in U$. 

\begin{lem} \mlabel{lem:5.7} 
If the function $\phi$ is dilatable, then 
the Hilbert space $\cH_\phi$ is the orthogonal direct sum of the subspaces 
$\cH_{\phi_n}$ and the series $\phi = \sum_n \phi_n$ converges uniformly on every ball %\break 
$B_r := \{a \in \cA \: p(a) < r\}$. 
\end{lem}

\begin{prf} Let $\cH := \cH_\phi$ and write $\pi$ for the corresponding representation 
of $S$ on~$\cH$ defined by $(\pi(s)f)(t) = f(ts)$. We write $\pi^1$ for the canonical 
extension of $\pi$ to $\ball(\cA^1,p^1)$ (Example~\ref{ex:3.14}). 
For $w := K_{s,v}$ and $t \in \cJ := (-1,1)$, we have 
\[\la \pi^1(t) w , w \ra = \la \pi^1(t) K_{s,v}, K_{s,v} \ra 
= \la K_{s,v}(st), v \ra = \la K_{st} K_s^* v, v \ra = \la \phi(sts^*)v,v\ra,\] 
and this function on $\cJ$ is analytic and positive definite, hence completely positive by 
Proposition~\ref{prop:5.7}. We conclude that $w$ has in $\cH$ an orthogonal expansion 
$w = \sum_n w_n$ with $\pi^1(t) w_n = t^n w_n$ for $t \in \cJ$. 
As the vectors of the form $K_{s,v}$, $s \in S$, $v \in V$, span a dense subspace of $\cH$, 
it follows that $\cH$ is the orthogonal direct sum of the $\cJ$-eigenspaces 
\[ \cH_n := \{ v \in \cH \: (\forall t \in \cJ)\ \pi(t) v = t^n v \}.\] 
Let $\pi_n$ denote the corresponding subrepresentation of $S$ on $\cH_n$. 
Accordingly, we write $\iota(v) = \sum_n \iota_n(v)$ with continuous linear maps 
$\iota_n \: V \to \cH_n$. This leads to the weakly convergent series 
\[ \phi(s) = \iota^* \pi(s) \iota
= \sum_n \iota_n^* \pi_n(s) \iota_n. \] 
We further have $\|\iota_n\| \leq \|\iota\|$ and, for 
$a \in \cA$ with $p(a)< r < 1$, 
\[ \|\pi_n(a)\|  
= r^n \|\pi_n(r^{-1}a)\| 
\leq r^n \|\pi(r^{-1}a)\| \leq r^n.\] 
Therefore the above series converges uniformly on every ball $B_r$, $r < 1$. 
Further, the homogeneity of the functions $\iota_n^* \pi_n(s) \iota_n$ implies that 
$\phi_n(s) = \iota_n^* \pi_n(s) \iota_n$ for every $s \in S$. 
This in turn implies that $\cH_n = \cH_{\phi_n}$.  
\end{prf}

\subsection{Homogeneous positive definite functions} 
\mlabel{subsec:6.3}

We now take a closer look at the representations corresponding to homogeneous positive definite functions 
$\phi$ of degree $n$ on $\ball(\cA,p)$, resp., $\cA$. In particular, we shall see that these functions 
are completely positive, and this will complete the proof of the result that bounded 
analytic positive definite functions are completely positive.

\begin{lem}\mlabel{lem:tens1}
If $(\cA_j, p_j)$, $j =1,2$, are seminormed involutive algebras, 
then their projective 
tensor product $(\cA_1 \otimes \cA_2, p_1 \otimes p_2)$ with 
\[ (p_1 \otimes p_2)(z) := \inf 
\Big\{ \sum_j p_1(x_j) p_2(y_j) \: z = \sum_j x_j \otimes y_j\Big\}\] 
is also a seminormed involutive algebra. 
\end{lem}

\begin{prf} It is well-known that $p_1 \otimes p_2$ defines a seminorm on the real associative 
$*$-algebra $\cB := \cA_1 \otimes \cA_2$ (\cite{Tr67}). 
For $z = \sum_j x_j \otimes y_j$ and 
$z' = \sum_j x_j' \otimes y_j'$, we have  
$zz' = \sum_{j,k} x_j x_k' \otimes y_j y_k'$. From 
\[ \sum_{j,k} p_1(x_j x_k') p_2(y_j y_k') 
\leq \sum_{j,k} p_1(x_j) p_2(y_j)  p_1(x_k') p_2(y_k')
= \sum_j  p_1(x_j) p_2(y_j)  \sum_k  p_1(x_k') p_2(y_k') \] 
we then derive that 
\[ (p_1 \otimes p_2)(zz') \leq (p_1 \otimes p_2)(z)(p_1 \otimes p_2)(z'),\]
i.e., that $p_1 \otimes q_2$ is submultiplicative. It is also easy to see that 
$p_1 \otimes p_2$ is involutive. 
\end{prf}

\begin{ex} The preceding lemma implies 
in particular, that we obtain for any 
seminormed involutive algebra $(\cA,p)$ a natural seminorm on the matrix algebra 
\break $M_n(\cA) \cong M_n(\R) \otimes_\R \cA$, where the norm on $M_n(\R)$ is the 
operator norm with respect to the euclidean norm on $\R^n$. 
\end{ex}

We consider the involutive algebra 
\[ S^n(\cA) := (\cA^{\otimes n})^{S_n} \] 
and endow it with the seminorm $p_n$ obtained by restricting the 
(projective) seminorm $p^{\otimes n}$ (Lemma~\ref{lem:tens1}). 
Then the natural homomorphism $\eta_{\cA}^{\otimes n} \: \cA^{\otimes n} \to C^*(\cA,p)^{\otimes n}$ 
intertwines the $S_n$-actions, so that it induces a morphism 
\[ S^n(\eta_{\cA})  \: S^n(\cA) \to S^n(C^*(\cA,p)).\] 

If $\cB$ and $\cC$ are $C^*$-algebras, we write 
$\cB \otimes \cC$ for the $C^*$-algebra obtained from the 
algebraic tensor product by completion with respect to the maximal 
$C^*$-norm, which leads to the 
{\it projective $C^*$-tensor product} (cf.\ \cite[Def.~IV.4.5]{Tak02}).

The following proposition generalizes \cite[Thm.3]{Ok66} that deals with 
the special case of Banach $*$-algebras with approximate identities.

\begin{prop} \mlabel{prop:tens} 
Let $(\cA_j, p_j)$, $j =1,2$, be two real seminormed involutive algebras. 
Then the map $\eta_{\cA_1} \otimes \eta_{\cA_2} \: \cA_1 \otimes \cA_2 \to 
C^*(\cA_1, p_1) \otimes  C^*(\cA_2, p_2)$ 
has the universal property of the enveloping $C^*$-algebra 
of the seminormed involutive algebra $(\cA_1 \otimes \cA_2, p_1 \otimes p_2)$. 
\end{prop}

\begin{prf} Let $\pi \: \cA_1 \otimes \cA_2 \to B(\cH)$ be a non-degenerate 
$p_1 \otimes p_2$-bounded representation. If the algebras $\cA_j$ 
are both unital, then $\pi_1(a_1) := \pi(a_1 \otimes \1)$ and 
$\pi_2(a_2) := \pi(\1 \otimes a_2)$ are pairwise commuting representations of 
$\cA_1$, resp., $\cA_2$ with 
\[ \pi(a_1 \otimes a_2) = \pi_1(a_1)\pi_2(a_2). \] 

If the algebras $\cA_j$ are not unital, we have to invoke multiplier 
techniques to obtain the representations $\pi_j$. 
Let $(\cA_j^1,p_j^1)$ denote the unital seminormed involutive 
algebra associated to $(\cA_j, p_j)$ 
(Definition~\ref{def:A1}). Then $\cA_1 \otimes \cA_2$ is an ideal 
of $\cA_1^1 \otimes \cA_2^1$, so that the latter algebra 
acts on $\cA_1 \otimes \cA_2$ by multipliers. 
For $m_j \in \cA_j^1$ and $z \in \cA_1 \otimes \cA_2$, we have  
\begin{align*}
(p_1 \otimes p_2)((m_1 \otimes m_2)z)
&= \inf \Big\{ \sum_k p_1(x_k) p_2(y_k) \: \sum_k x_k \otimes y_k = (m_1 \otimes m_2)z \Big\} \\
&\leq \inf \Big\{ \sum_k p_1(m_1 x_k) p_2(m_2 y_k) \: \sum_k x_k \otimes y_k = z \Big\} \\
&\leq p_1^1(m_1) p_2^1(m_2) \inf \Big\{ \sum_k p_1(x_k) p_2(y_k) \: \sum_k x_k \otimes y_k 
= z\Big\} \\
& = p_1^1(m_1) p_2^1(m_2) (p_1 \otimes p_2)(z).
\end{align*}
We conclude that, 
for $S_j := \oline\ball(\cA_1^1, p_1^1)$ and $S := 
 \ball(\cA_1 \otimes \cA_2, p_1 \otimes p_2)$, we have 
\[ (S_1 \times S_2) S \subeq S.\] 
Applying Lemma~\ref{lem:ideal-lin-ext} to the ideal 
\[ \cA := \cA_1 \otimes \cA_2 
\trile \cB := (\cA_1 \otimes \cA_2) + \cA_1^1 \otimes \1, \] 
we obtain a unique linear representation $\pi_1^1 \: \cA_1^1 \to B(\cH)$, 
bounded on $S_1$ and satisfying 
\[ \pi_1^1(a)\pi(a_1 \otimes a_2) =  \pi(aa_1 \otimes a_2) 
\quad \mbox{ for } \quad a, a_1 \in \cA_1, a_2 \in \cA_2.\] 
We likewise obtain a unique linear representation $\pi_2^1 \: \cA_2^1 \to B(\cH)$, 
bounded on $S_2$ and satisfying 
\[ \pi_2^1(a)\pi(a_1 \otimes a_2) =  \pi(a_1 \otimes aa_2) 
\quad \mbox{ for } \quad a_1 \in \cA_1, a, a_2 \in \cA_2.\] 
We put $\pi_j := \pi_j^1\res_{\cA_j}$. For $a_1, a_1' \in \cA$ and 
$a_2, a_2' \in \cA_2$, the relation 
\[ \pi_1(a_1) \pi_2(a_2) \pi(a_1' \otimes a_2') 
= \pi(a_1 a_1' \otimes a_2 a_2') 
= \pi(a_1 \otimes a_2) \pi(a_1' \otimes a_2') \] 
now implies that 
\[ \pi(a_1 \otimes a_2) =  \pi_1(a_1) \pi_2(a_2).\] 
Since $\pi_j^1$ is bounded on $S_j$, it follows that the 
representation $\pi_j$ of $\cA_j$ is bounded on $\ball(\cA_j, p_j)$, i.e., 
$\|\pi_j(a_j)\| \leq p_j(a_j)$ for $a_j \in \cA_j$, $j =1,2$. 
Hence there exist pairwise commuting representations 
$\hat\pi_j \: C^*(\cA_j,p_j) \to B(\cH)$ with $\hat\pi_j \circ \eta_{\cA_j} = \pi_j$. Now 
\[ \hat\pi(a_1 \otimes  a_2) := 
\hat\pi_1(a_1)\hat\pi_2(a_2) \] 
defines a representation $\hat\pi \: C^*(\cA_1,p_1) \otimes C^*(\cA_2, p_2) 
\to B(\cH)$ with $\hat\pi \circ (\eta_{\cA} \otimes \eta_{\cA_2})=\pi$ 
(see \cite[Prop.~IV.4.7]{Tak02}). 
\end{prf}

Iterating the preceding lemma, we obtain in particular: 
\begin{cor} \mlabel{cor:ten} 
The map $\eta_{\cA}^{\otimes n}\: \cA^{\otimes n} \to C^*(\cA, p)^{\otimes n}$ 
is universal for the seminormed involutive algebra $(\cA^{\otimes n}, p^{\otimes n})$. 
\end{cor}

\begin{lem} \mlabel{lem:fixproj} Let $\Gamma$ be a finite group acting on the seminormed 
involutive algebra $(\cA,p)$. Then the fixed point projection 
$p \: \cA \to \cA^\Gamma$ 
is completely positive. 
\end{lem}

\begin{prf} Since $M_n(p) \: M_n(\cA) \to M_n(\cA^\Gamma)$ is the fixed point 
projection for the canonical action of $\Gamma$ on $M_n(\cA)$, it suffices to 
show that $p$ is positive. This follows from the positivity of 
\[ p(a^*a) 
= \frac{1}{|\Gamma|} \sum_{\gamma \in \Gamma} \gamma(a^*a)
= \frac{1}{|\Gamma|} \sum_{\gamma \in \Gamma} \gamma(a)^*\gamma(a).\qedhere\] 
\end{prf}

\begin{prop} \mlabel{prop:sn-univ} 
The homomorphism $S^n(\eta_\cA) \: S^n(\cA) \to S^n(C^*(\cA,p))$ has the 
universal property of $C^*(S^n(\cA), p_n)$. 
\end{prop}

\begin{prf} {\bf Case 1:} First we assume that $\cA$ is unital. 
Let $\pi \: S^n(\cA) \to B(\cH)$ be a $p_n$-bounded (unital) representation. 
The $S_n$-equivariant projection 
$q \: \cA^{\otimes n} \to S^n(\cA)$ is completely positive  
(Lemma~\ref{lem:fixproj}) and since the seminorm $p^{\otimes n}$ is $S_n$-invariant, 
$q$ is a $p^{\otimes n}$-contraction. Therefore 
$\pi \circ q \: \cA^{\otimes n} \to B(\cH)$ is $p^{\otimes n}$-bounded and 
completely positive, hence of the form 
$\iota^* \rho(\cdot) \iota$ for a $p^{\otimes n}$-bounded representation 
$\rho$ of $\cA^{\otimes n}$  (cf.~Theorem~\ref{thm:Arv2.2}). 
In view of Proposition~\ref{cp-univ} and Corollary~\ref{cor:ten}, 
this implies the existence of a 
completely positive map 
\[ \Phi \:  C^*(\cA,p)^{\otimes n} \to B(\cH) \quad \mbox{ with } \quad 
\Phi(\eta_\cA(a_1) \otimes \cdots \otimes \eta_\cA(a_n)) 
= \pi(q(a_1 \otimes \cdots \otimes a_n)).\] 
Then 
$\Phi(\eta_\cA(a) \otimes \cdots \otimes \eta_\cA(a)) 
= \pi(a \otimes \cdots \otimes a))$ implies that 
\[ \Phi \circ S^n(\eta_\cA) = \pi.\] 
As $\pi$ is a representation and the image of $S^n(\eta_\cA)$ spans a dense subspace in 
$S^n(C^*(\cA,p))$, 
$\Phi$ is multiplicative on $S^n(C^*(\cA,p))$, so that 
$\Phi\res_{S^n(\cA)}$ is a representation. 

{\bf Case 2:} If $\cA$ is not unital, we assume that $\pi$ is non-degenerate. 
We first observe that 
$S^n(\cA) \trile S^n(\cA^1)$ is an ideal. 
Let 
\[ S := \ball(S^n(\cA), p_n) \quad \mbox{ and } \quad 
S^1 := \{ b \in S^n(\cA^1) \: b S \cup S b \subeq S \}.\] 
Then $S^1$ contains all elements of the form 
$a^{\otimes n}$, where $a \in \ball(\cA^1, p^1)$, and these elements 
span $S^n(\cA^1)$. Therefore Lemma~\ref{lem:ideal-lin-ext} implies that 
$\pi$ extends to a linear unital representation 
$\hat\pi \: S^n(\cA^1) \to B(\cH)$ which is bounded on $S^1$. 

We claim that $\hat\pi$ is bounded on $\ball(S^n(\cA^1), p^1_n)$. 
To this end, we first observe that, for $a \in \cA$, $m \in \cA^1$, we have
$p(ma) \leq p^1(m) p(a)$ because $p^1$ is submultiplicative and extends $p$. 
This leads to 
\[ p^{\otimes n}((m_1 \otimes \cdots \otimes m_n)z) 
\leq p^1(m_1) \cdots p^1(m_n) p^{\otimes n}(z) \quad \mbox{ for } \quad 
m_j \in \cA^1, z \in \cA^{\otimes n},\] 
and this in turn to 
\[ p^{\otimes n}(wz) \leq (p^1)^{\otimes n}(w) p^{\otimes n}(z) 
\quad \mbox{ for } \quad 
w \in (\cA^1)^{\otimes n}, z \in \cA^{\otimes n}.\] 
This specializes to 
\[ p_n(wz) \leq p_n^1(w) p_n(z) \quad \mbox{ for } \quad 
w \in S^n(\cA^1), z \in S^n(\cA).\] 
From this estimate we derive that 
$\ball(S^n(\cA^1),p^1_n) \subeq S^1$, so that 
$\hat\pi$ is bounded on the semigroup $\ball(S^n(\cA^1),p^1_n)$, 
which leads to $\|\hat\pi(w)\| \leq p^1_n(w)$ for 
$w \in S^n(\cA^1)$. Hence Case $1$ shows that there exists a 
unital representation 
$\beta \: S^n(C^*(\cA^1,p)) \to B(\cH)$ with 
\[ \beta \circ \eta_{S^n(\cA^1)} = \hat\pi.\] 
Composing with the natural map 
$S^n(C^*(\cA,p)) \to S^n(C^*(\cA,p)^1) \cong S^n(C^*(\cA^1,p^1))$ 
(cf.\ Remark~\ref{rem:unitization}), we obtain a representation 
$\gamma \: S^n(C^*(\cA,p)) \to B(\cH)$ with 
\[ \gamma(S^n(\eta_\cA)(a^{\otimes n})) 
= \gamma(\eta_\cA(a)^{\otimes n}) 
=  \beta(\eta_{\cA^1}(a)^{\otimes n}) 
=  \beta(\eta_{S^n(\cA^1)}(a^{\otimes n})) = \hat\pi(a^{\otimes n}) 
= \pi(a^{\otimes n}).\] 
This implies that the two representations 
$\pi$ and $\gamma \circ S^n(\eta_\cA)$ of $S^n(\cA)$ coincide, 
and thus $S^n(\eta_\cA) \: S^n(\cA) \to S^n(C^*(\cA,p))$ 
has the universal property of the enveloping $C^*$-algebra of the seminormed 
involutive algebra $(S^n(\cA), p_n)$. 
\end{prf}

The following automatic continuity lemma 
is a version of Corollary~\ref{cp-univ} for polynomials of higher degree. 
We recall that the dilatability assumption is always satisfied if $\cA$ is 
unital (Theorem~\ref{thm:Arv2.2}). 

\begin{lem}\mlabel{Arv3.20}
Let $(\cA,p)$ be a real seminormed involutive algebra, 
$n \in \N_0$, and $V$ be a complex Hilbert space. 
For every $p$-bounded dilatable positive definite homogeneous polynomial 
$\varphi\colon \cA\to B(V)$ of degree $n$,  
there exists a unique complex linear completely positive map 
\[ \Phi_n\colon S^n(C^*(\cA,p))\to B(V) \quad \mbox{ with } \quad 
\varphi(a)=\Phi_n(\eta_{\cA}(a)^{\otimes n})
\quad \mbox{ for } \quad a\in \cA.\] 
In particular, $\phi$ is completely positive. 
\end{lem}

\begin{prf} Let $(\pi, \cH, \iota)$ be a dilation of $\phi\res_{\ball(\cA,p)}$. 
The assumption that $\phi$ is a homogeneous polynomial of degree $n$ implies that 
all elements of $\cH \cong \cH_\phi$ have this property, so that 
$\pi(r a) := r^n \pi(a)$ yields an extension to a representation 
of the involutive semigroup $(\cA,\cdot)$. % with $\pi(r) = r^n \1$ for $r \in \R$. 
The map $\pi \: \cA \to B(\cH)$ is a 
polynomial of degree $n$ bounded on $\ball(\cA,p)$. Hence 
$\pi$ defines a unique linear representation 
$\tilde\pi \: S^n(\cA) \to B(\cH)$ bounded on the subsemigroup 
\[ S := \{ a^{\otimes n}\: a \in \ball(\cA,p)\} \subeq \ball(S^n(\cA),p_n).\] 
This implies that the corresponding symmetric $n$-linear map 
$\beta \: \cA^n \to B(\cH)$ is continuous on the diagonal with respect to the 
$p$-topology, and polarization shows that $\beta$ is continuous with respect to the 
$p$-topology on $\cA^n$. This in turn shows that the corresponding 
linear map $\tilde\beta \: \cA^{\otimes n} \to B(\cH)$ is $p^{\otimes n}$-continuous, 
and therefore $\tilde\beta\res_{S^n(\cA)} = \tilde\pi$ is 
$p_n$-continuous, i.e., bounded on $\ball(S^n(\cA),p_n)$, and therefore 
$\|\tilde\pi(z)\| \leq 1$ for $z \in \ball(S^n(\cA),p_n)$ 
(cf.\ Remark~\ref{rem:cstar}). 

Now Proposition~\ref{prop:sn-univ} provides a 
representation $\rho \: S^n(C^*(\cA,p)) \to B(\cH)$ with 
$\rho \circ S^n(\eta_\cA) = \tilde\pi$. We thus obtain 
for $a \in \cA$: 
\[ \phi(a) 
= \iota^* \pi(a) \iota 
= \iota^* \tilde\pi(a^{\otimes n}) \iota 
= \iota^* \rho(\eta_\cA(a)^{\otimes n}) \iota.\] 
Hence $\phi = \Phi_n(\eta_{\cA}(a)^{\otimes n})$ 
holds for $\Phi_n(b) := \iota^* \rho(b) \iota$, which is obviously completely positive, 
hence also positive definite. 

For $\Phi(b) := \Phi_n(b^{\otimes n})$ we now have 
$\Phi \circ \eta_\cA = \phi$. Therefore the complete positivity 
of $\phi$ follows from the complete positivity of $\eta_\cA$, 
the complete positivity of the map \break $C^*(\cA,p) \to S^n(C^*(\cA,p)), b \mapsto 
b^{\otimes n}$ (\cite[Lemma~3.5]{Arv87}) and the complete positivity of~$\Phi_n$. 
\end{prf}

\begin{thm} \mlabel{thm:5.13} 
If $(\cA,p)$ is a seminormed algebra, $V$ a complex Hilbert space and \break 
$\phi \: \ball(\cA,p) \to B(V)$ a dilatable bounded  positive definite function which 
is analytic with respect to some locally convex topology for which $p$ is continuous, 
then $\phi$ is completely positive. 
\end{thm}

\begin{prf} First, Lemma~\ref{lem:5.7} 
asserts that the Taylor series $\phi = \sum_n \phi_n$ converges uniformly on 
all balls $B_r$, $r < 1$. Next, Lemma~\ref{Arv3.20} implies that the functions 
$\phi_n$ are completely positive, so that $\phi$ is completely positive as well. 
\end{prf}

\section{Proof of the main theorem} 
\mlabel{sec:6} 

In this section we finally prove Theorem~\ref{thm:intro}. 

\begin{defn}\label{exp_alg}
For any $C^*$-algebra $\cB$ and $k \in \N_0$, we consider the 
$C^*$-algebra 
\[ S^n(\cB) := (\cB^{\otimes n})^{S_n} \] 
generated by the elements of the form $b^{\otimes n}$, $b \in \cB$, in the 
$C^*$-algebra $\cB^{\otimes n}$. Here the tensor product means the 
{\it projective $C^*$-tensor product}, resp., the completion of the tensor 
product with respect to the maximal $C^*$-norm 
(cf.\ \cite[Def.~IV.4.5]{Tak02}). 
The above $S^n(\cB)$ is the subalgebra of fixed points for the 
natural action of the symmetric group $S_n$ on $\cB^{\otimes n}$. 
For $n = 0$, we put $S^0(\cB):=\C$. 
The \emph{exponential $C^*$-algebra} of $\cB$ is 
\begin{equation}\label{Arv_above3.22}
e^{\cB}:=\Bigl\{(b_k)_{k\ge 0}\in\prod\limits_{k\ge 0}S^k(\cB)\mid\lim\limits_{k\to\infty}\Vert b_k\Vert=0\Bigr\}.\end{equation}
We then obtain a natural map 
\[ \Gamma_\cB \colon \ball(\cB)\to \oline{\ball}(e^{\cB}), 
\quad \Gamma_\cB(b)=\sum\limits_{k\ge 0} b^{\otimes k}.\]
We will write simply $\Gamma_\cB=\Gamma$ whenever the $C^*$-algebra $\cB$ is understood from the 
context. 
\end{defn}

\begin{lem}\label{CP}
For every $C^*$-algebra $\cB$, the map $\Gamma_\cB\colon \ball(\cB)\to C(e^{\cB})$ 
is a completely positive holomorphic homomorphism of $*$-semi\-groups.  
\end{lem}

\begin{prf}
For every $k \in \N_0$, the map $\Gamma_k\vert_{\ball(\cB)}\colon \ball(\cB)\to S^k(\cB)$ 
is completely positive and homogeneous of degree $k$ by \cite[Lemma 3.5]{Arv87}. 
Since their sum converges uniformly on every ball of radius $r< 1$, 
the map $\Gamma_\cB$ is also holomorphic and  completely positive. 
\end{prf}

We now consider the natural map 
\[ \Gamma = \Gamma_\cA \: \ball(\cA,p)\to \ball(e^{C^*(\cA,p)}), 
\quad \Gamma(a)=\sum\limits_{k\ge 0} \eta_\cA(a)^{\otimes k}.\] 

\begin{lem} \mlabel{lem:1.7} 
The map $\Gamma \: \ball(\cA,p) \to e^{C^*(\cA,p)}$ 
is completely positive and 
analytic in the $p$-topology on $\cA$. 
\end{lem}

\begin{prf} Let $\cB := C^*(\cA,p)$. 
Since the map $\eta_\cA \: \ball(\cA,p) \to \ball(\cB)$ 
is continuous linear, the analyticity of $\Gamma$ 
follows from the holomorphy of $\Gamma_\cB$ 
(Lemma~\ref{CP}). 

To see that $\Gamma$ is completely positive, we % of type (S), we 
have to show that, for $A = (a_{ij}) \geq 0$ in $M_n(\cA)$ and 
$p(a_{ij})< 1$, the matrix $(\Gamma(a_{ij})) \in M_n(e^\cB)$ is positive. 
This follows from the complete positivity of the linear $*$-homomorphism 
$\eta_\cA \: \cA \into \cB$ (Remark~\ref{rem:1.5}(b)) and Lemma~\ref{CP}. 
\end{prf}

The following theorem is our main theorem stated as 
Theorem~\ref{thm:intro} in the introduction.

\begin{thm}\label{thm:6.3} 
Let $(\cA,p)$ be a real seminormed involutive algebra, $V$ a complex Hilbert space and 
$\phi \: \ball(\cA,p) \to B(V)$ be a bounded function. Then the following are equivalent: 
\begin{itemize}
\item[\rm(i)] $\phi$ is completely positive and dilatable. 
\item[\rm(ii)] $\phi$ is dilatable, positive definite and analytic 
with respect to some locally convex topology for which $p$ is continuous. 
\item[\rm(iii)] There exists a linear completely positive map 
$\Phi \: e^{C^*(\cA,p)}\to B(V)$ with $\Phi \circ \Gamma = \phi$. 
\end{itemize}
\end{thm}

\begin{prf} (i) $\Rarrow$ (ii) follows from Theorem~\ref{thm:4.7}. 

(ii) $\Rarrow$ (i) follows from Theorem~\ref{thm:5.13}. 

(iii) $\Rarrow$ (i): Let $\cB := C^*(\cA,p)$ and  
$\Phi\colon e^{\cB}\to B(V)$ be any completely positive linear map. 
Since a composition of completely positive maps is completely positive 
(Remark~\ref{rem:1.5}), 
it follows from Lemma~\ref{lem:1.7} that $\Phi\circ {\Gamma}$ is an analytic 
completely positive map. Further, the dilatability of $\phi$ follows from 
the dilatability of $\Phi$ (cf.\ Proposition~\ref{cp-univ}). 

(i) $\Rarrow$ (iii): First, Theorem~\ref{thm:4.7} implies the existence 
of an expansion $\phi = \sum_{n = 0}^\infty \phi_n$, where 
each $\phi_n$ is a homogeneous completely positive polynomial 
with $\|\phi_n\|_\infty \leq \|\phi\|_\infty$. 
Lemma~\ref{Arv3.20} now provides linear positive maps 
$\Phi_n \:  S^n(C^*(\cA,p)) \to B(V)$ with 
$\phi_n(a) = \Phi_n(\eta_\cA(a)^{\otimes n})$ for 
$a \in \ball(\cA,p)$. For the canonical extensions to $\cA^1$, we further have 
$\Phi_n^1(\1^{\otimes n}) = \phi_n^1(\1) \leq \phi^1(\1)$, so that 
the sequence $(\|\Phi_n^1\|)_{n \in \N_0}$ is bounded. This implies that 
$\Phi(x) := \sum_n \Phi_n(x_n)$ defines a continuous linear function 
$e^{C^*(\cA,p)} \to B(V)$ with the required properties. 
\end{prf}

For unital algebras, we obtain the  following simplification:

\begin{cor}\label{cor:Arv3.1}
Let $(\cA,p)$ be a real unital 
seminormed involutive algebra, $V$ a complex Hilbert space and 
$\phi \: \ball(\cA,p) \to B(V)$ be a bounded function. Then the following are equivalent: 
\begin{itemize}
\item[\rm(i)] $\phi$ is completely positive. 
\item[\rm(ii)] $\phi$ is positive definite and analytic 
with respect to some  locally convex topology for which $p$ is continuous. 
\item[\rm(iii)] There exists a linear completely positive map 
$\Phi \: e^{C^*(\cA,p)}\to B(V)$ with $\Phi \circ \Gamma = \phi$. 
\end{itemize}
\end{cor}

\begin{prf} If $\cA$ is unital, then every bounded positive definite 
function on $\ball(\cA,p)$ is dilatable by Theorem~\ref{thm:Arv2.2}, 
so that the assertion follows from Theorem~\ref{thm:6.3}. 
\end{prf}

Recall from Example~\ref{ex:non-dil} that the dilatability assumptions cannot 
be removed from Theorem~\ref{thm:6.3}, which is the same as Theorem~\ref{thm:intro}. 

\begin{cor} The assignment $\Phi\mapsto\Phi\circ \Gamma$ 
is a bijection between the completely positive linear $B(V)$-valued maps 
on the $C^*$-algebra $e^{C^*(\cA,p)}$ 
and the bounded dilatable completely positive $B(V)$-valued 
functions on the $*$-semigroup~$\ball(\cA,p)$. 
\end{cor}

\begin{prf} That $\Phi \circ \Gamma$ determines $\Phi$ uniquely 
follows from the fact that 
the complex linear span of $\Gamma(\ball(\cA,p))$ is dense in 
$e^{C^*(\cA,p)}$ (\cite[Lemma 3.23]{Arv87}).  
\end{prf}

\begin{rem} If the real seminormed involutive algebra 
$(\cA,p) = (\cA,\|\cdot\|)$  is a unital $C^*$-algebra, 
then $C^*(\cA,p) \cong \cA_{\C}\cong \cA\oplus\oline\cA$  
and the completely positive linear maps on 
$\Phi \: e^{\cA\oplus\oline\cA}\cong e^{\cA}\otimes e^{\oline\cA}\to \C$ normalized 
by $\|\Phi\| = 1$ are precisely the states of $e^{\cA}\otimes e^{\oline\cA}$. 
This special case of the above Theorem~\ref{thm:6.3} is \cite[Th. 3.1]{Arv87}. 
\end{rem}

The following theorem is an important consequence of Theorem~\ref{thm:6.3} 
for the representation theory of $\ball(\cA,p)$. 

\begin{thm}\label{final0}
Let $(\cA,p)$ be a real seminormed involutive algebra and $\cH$ a complex 
Hilbert space.  Then the correspondence 
$\hat\pi\mapsto \hat\pi\circ \Gamma$ 
defines a commutant-preserving bijection between the $*$-representations of the $C^*$-algebra 
$e^{C^*(\cA,p)}$ on $\cH$ and the bounded real-analytic $*$-representations 
of the $*$-semigroup $\ball(\cA,p)$ on~$\cH$. 
\end{thm}

\begin{prf}
Every $*$-representation of a $*$-semigroup is a dilatable positive definite function, 
hence we can use Theorem~\ref{thm:6.3} to see that 
a bounded real-analytic unital $*$-representation \break 
$\pi\colon \ball(\cA)\to B(\cH)$ 
is completely positive and that 
there exists a unique completely positive linear map 
$\Phi\colon e^{C^*(\cA,p)}\to B(\cH)$ with 
$\pi=\Phi\circ {\Gamma}.$ 
It remains to prove that, for all $x,y\in e^{C^*(\cA,p)}$, we have $\Phi(xy)=\Phi(x)\Phi(y)$. 
Since $\Phi$ is positive, it is also continuous even though it is defined on the 
non-unital $C^*$-algebra~$e^{C^*(\cA,p)}$ (see \cite[Ch. 2]{Pau02}). 
As $\pi=\Phi\circ\Gamma$ and both $\pi$ and ${\Gamma}$ are semigroup homomorphisms, 
$\Phi(xy)=\Phi(x)\Phi(y)$ holds for $x,y\in\Gamma(\ball(\cA,p))$, 
and since $\Gamma(\ball(\cA,p))$ spans a dense linear subspace of $e^{C^*(\cA,p)}$ 
and $\Phi$ is continuous, it is also multiplicative, 
hence a $*$-representation of the $C^*$-algebra $e^{C^*(\cA,p)}$. 

Conversely, if $\Phi$ is a $*$-homomorphism of $C^*$-algebras, 
then $\Phi\circ \Gamma$ 
is a bounded real-analytic $*$-representation 
of the $*$-semigroup $\ball(\cA,p)$ on $\cH$.
\end{prf}

\begin{rem} \mlabel{rem:6.8}
(a) Theorem~\ref{final0} asserts that the bounded real-analytic representation 
theory of the involutive semigroup $\ball(\cA,p)$ is faithfully represented 
by the $C^*$-algebra $e^{C^*(\cA,p)}$, so this is a host algebra in the 
sense of \cite{Gr05}. Since the real linear map 
$\eta_\cA \: \cA \to C^*(\cA,p)$ extends to a complex linear map 
$\cA_\C \to C^*(\cA,p)$, it follows in particular 
that every bounded analytic completely positive function 
$\phi \: \ball(\cA,p) \to B(V)$ extends to a holomorphic function 
on the open semigroup $\ball(\cA_\C, p_\C)$ for 
$p_\C(a+ib) := p(a) + p(b)$ (cf.\ Definition~\ref{def:complexif}). 
Therefore $e^{C^*(\cA,p)}$ is also a host algebra for the bounded holomorphic 
representations theory of the complex semigroup $\ball(\cA_\C,q)$ in 
the sense of \cite{Ne08}. 

(b) Since $e^{C^*(\cA,p)}$ is a $c_0$-direct sum of the ideals 
$S^N(C^*(\cA,p))$, all its non-degenerate 
representations decompose as direct sums of representations 
on which exactly one of these ideals acts non-trivially. Therefore the main 
point in applying Theorem~\ref{final0} is to identify the representations 
of $S^N(C^*(\cA,p)) \cong C^*(S^N(\cA),p_N)$ in terms of the (linear) representations 
of the $C^*$-algebra $C^*(\cA,p)$, resp., representations of the algebra~$\cA$ which are 
bounded on $\ball(\cA,p)$. 

(c) For every $N \in \N$, we see that 
$C^*$-algebra $S^N(C^*(\cA,p))$ is a host algebra for the 
class of the bounded $N$-homogeneous polynomial representations of $\ball(\cA,p)$, resp., 
the $N$-homoge\-ne\-ous polynomial representations of $\cA$ bounded on $\ball(\cA,p)$. 
\end{rem}

\begin{cor}\label{final1}
Let $\cA$ be a $C^*$-algebra, considered as a real seminormed 
algebra, and $\cH$ be a complex Hilbert space.   
Every bounded real-analytic $*$-representation $\pi\colon \ball(\cA)\to B(\cH)$
is completely positive and the correspondence 
$\Phi\mapsto\Phi\circ(\Gamma_\cA\otimes\Gamma_{\oline\cA})$ 
is a bijection between the $*$-representations of the $C^*$-algebra 
$e^{\cA}\otimes e^{\oline\cA}$ on $\cH$ 
and the bounded real-analytic $*$-representations of the $*$-semigroup 
$\ball(\cA)$ on~$\cH$. 
\end{cor}

\begin{prf}
This is the special case of Theorem~\ref{final0}, where $\cA$ is a unital $C^*$-algebra. 
\end{prf}

\begin{ex} (cf.\ \cite[Ex.~2.1]{HN87}) 
We consider the non-unital $C^*$-algebra $\cA = c_0(\N_0,\C)$. 
Then 
\[ \phi(z) := \sum_{n = 0}^\infty z_n^n \] 
is a completely positive $\C$-valued function on $\cA$ which is bounded on 
the ball $B_s(0)$ of radius $s$ if and only if $s < 1$. In particular, it is unbounded 
on $\ball(\cA)$. For $r > 0$ we put $\phi_r(z) := \phi(rz)$. Then 
$\phi_r$ is bounded on $\ball(\cA)$ if $r < 1$, hence dilatable because 
$\cA$ is a $C^*$-algebra (Proposition~\ref{cp-univ}). 

In view of 
\[ \phi(zw^*) := \sum_{n = 0}^\infty z_n^n \oline{w_n}^n,\] 
the functions $\phi_n(z) = z_n^n$ form an orthonormal basis of the 
Hilbert space $\cH_\phi$, which implies in particlar that 
$\phi \not\in \cH_\phi$. For $r < 1$, we have $\|\phi_n\| = r^{-n/2}$, so that 
$\sum_n  r^{2n} \| \chi_n\|^2 = \sum_n  r^n < \infty$ implies that 
$\phi_r \in \cH_{\phi_r}$ holds for $r < 1$ (cf.\ Proposition~\ref{prop:exis-dil}). 
\end{ex}

\section{Relations to Arveson's c.p.\ concept} 
\mlabel{sec:7}

Let $\cA$ and $\cB$ be $C^*$-algebras. 
We say that a map $\phi \: \ball(\cA) \to \cB$ 
is {\it completely positive  of type (W)} 
if every positive element $A = (a_{ij}) \in M_n(\cA)$ with $\|A\| < 1$ 
is mapped to a positive element of $\cB$. 
This is the concept of complete positivity used in \cite{Arv87} for functions 
on the open unit ball of unital a $C^*$-algebra. 

\begin{rem} \mlabel{rem:w1} (a) If $\vphi \: \ball(\cA) \to \cB$ 
is homogeneous of degree~$\alpha\in\R$ as in Lemma~\ref{lem:1.5}, 
then complete positivity of type (W) implies complete positivity of its extension 
$\hat\phi \: \cA \to \cB$. In fact, for any $A = (a_{ij}) \in M_n(\cA)_+$, 
there exists an $r > 0$ 
with $r \|a_{ij}\| < 1$ for all $i,j$ and 
\[ \|(r a_{ij})\|= r\|(a_{ij})\| < 1.\] 
Then $(\hat\vphi(a_{ij})) = r^{-\alpha}(\vphi(r a_{ij})) \in M_n(\cB)_+$
implies that $\hat\vphi \: \cA \to \cB$ is completely positive. 

(b) We will show in Example~\ref{16jun-ex} that complete positivity of type (W) does not imply positive definiteness. 
This is related to the fact that the matrix $A$ used in Lemma~\ref{lem:3.2} need not be a contraction. 
For example the matrix 
\[ A := \pmat{r^2 & rs \\ rs & s^2}, \quad 0 < r,s < 1 \] 
satisfies $\|A\| = r^2 + s^2$.  
\end{rem} 

The next proposition shows that, for bounded positive definite functions on $\ball(\cA,p)$, 
complete positivity of type (W) implies complete positivity in the sense of Definition~\ref{defCP}(b). 

\begin{prop}\mlabel{prop:3.x} 
Let 
$\cA$ be a $C^*$-algebra 
and  
$V$  be a complex Hilbert space.  
If a bounded positive definite function $\varphi\colon \ball(\cA)\to B(V)$ is completely positive of type {\rm(W)}, 
then it is also completely positive. 
\end{prop}

\begin{prf} Let $\varphi\colon \ball(\cA)\to B(V)$ be completely positive of type {\rm(W)}. 
It follows by Remark~\ref{rem:w1}(c) that the corresponding GNS representation 
$(\pi_\vphi,\cH_\vphi)$ is completely positive of type (W). Since this representation is non-degenerate, 
it is a direct sum of cyclic representation $(\pi_j, \cH_j)_{j \in J}$ generated by 
cyclic vectors $v_j \in \cH_j \subeq \cH_\vphi$ (\cite{Ne00}). 
Then the functions 
\[ \vphi_j(s) 
:= \la \pi_j(s)v_j, v_j \ra 
= \la \pi(s)v_j, v_j \ra \] 
are also completely positive of type (W) since the linear functional 
$B(\cH_\vphi) \to \C, A \mapsto \la Av_j, v_j\ra$ is 
positive definite, hence completely positive (Remark~\ref{rem:1.5}(c)). 
Now Remark~\ref{rem:w1}(d) implies that the functions $\vphi_j$ are completely 
positive, so that the representations $\pi_j$ are completely positive 
by Proposition~\ref{prop:4.1}. This in turn implies that $\pi_\vphi$ is completely 
positive, which, again by Proposition~\ref{prop:4.1}, implies that 
$\vphi$ is completely positive. 
\end{prf}

We will show by a simple example that complete positivity of type (W) does not imply 
positive definiteness, hence also not complete positivity (see Example~\ref{16jun-ex}). 
The following simple observation is the key to that example. 
%Example~\ref{16jun-ex}. 

\begin{lem}\mlabel{16jun}
Let $V$ be any Hilbert space and $A\in B(V)$. 
If $0\le A\le \1$ and $u,v\in V$ with $\Vert u\Vert=\Vert v\Vert=1$ and $\la u,v\ra=0$, 
then $\vert\la Au,v\ra\vert\le1/2$. 
\end{lem}

\begin{prf}
Since $A\ge 0$, by the Cauchy--Schwartz inequality we obtain 
$$\vert\la Au,v\ra\vert^2\le \la Au,u\ra\la Av,v\ra.$$ 
Similarly, since $\1-A\ge 0$, we have $\vert\la (\1-A)u,v\ra\vert^2\le \la (\1-A)u,u\ra\la (\1-A)v,v\ra$,  
and then, by using the hypothesis on $u$ and $v$, we obtain 
$$\vert\la Au,v\ra\vert^2\le (1-\la Au,u\ra)(1-\la Av,v\ra).$$
By multiplying the above displayed inequalities and using the fact that $t(1-t)\le 1/4$ if $0\le t\le 1$, 
we obtain $\vert\la Au,v\ra\vert^4\le1/16$, which leads to the asserted inequality. 
\end{prf}

\begin{ex}\label{16jun-ex}
Let $\varphi\colon{\mathbb D}\to\C$ satisfying $\varphi(z)=z$ if $\vert z\vert\le 1/2$ or $z\in{\mathbb D}\cap[0,\infty)$. 
Then the following assertions hold:
\begin{enumerate}
\item For every integer $n\ge 1$ and $0\le A=(a_{ij})\in M_n(\C)$ with $\Vert A\Vert<1$ 
we have $0\le(\varphi(a_{ij}))\in M_n(\C)$. 
\item If there exists $z_0\in{\mathbb D}$ with $\varphi(\overline{z_0})\ne\overline{\varphi(z_0)}$, 
then there exist $a_1,a_2\in{\mathbb D}$ for which the matrix $(\varphi(a_i\overline{a_j}))_{1\le i,j\le 2}\in M_2(\C)$ 
fails to be nonnegative. 
\end{enumerate}

For the first assertion, if $0\le A=(a_{ij})\in M_n(\C)$ with $\Vert A\Vert<1$, then $0\le A\le\1$, 
hence by Lemma~\ref{16jun} we obtain $\vert a_{ij}\vert\le1/2$ if $i\ne j$. 
Moreover, $0\le a_{jj}\le 1$. 
Consequently $\varphi(a_{ij})=a_{ij}$ for all $i,j\in\{1,\dots,n\}$, 
and then trivially $0\le(\varphi(a_{ij}))\in M_n(\C)$. 

For the second assertion, select any $a_1,a_2\in{\mathbb D}$ with $a_1\overline{a_2}=z_0$. 
Since $\varphi(\overline{z_0})\ne\overline{\varphi(z_0)}$, it follows that the matrix 
$(\varphi(a_i\overline{a_j}))\in M_2(\C)$ 
fails to be hermitian, hence in particular it cannot be nonnegative. 
\end{ex}

%\framebox{add a comment on \cite{Arv87}}
The above example shows that \cite[Lemma 5.3]{Arv87} actually needs a stronger 
notion of complete positivity,  
which should ensure that every completely positive function in the corresponding sense 
is positive definite, so that \cite[Th. 2.2]{Arv87} is applicable in that setting. 
As we have seen in the present paper, the stronger notion of complete positivity from 
Definition~\ref{defCP}(b) 
not only fills that gap from \cite{Arv87} but also is flexible enough to allow 
one to develop dilation theory for maps defined on the unit ball of any real 
seminormed involutive algebra. 
As explained in the introduction, 
this level of generality is necessary for the applications to mapping algebras 
$\cA=C^\infty(X,\cB)$ for any $C^*$-algebra $\cB$, where $X$ is a smooth manifold 
(see \cite{BN15} for details).

\appendix
\section{Analytic functions on infinite dimensional spaces} 

In this appendix, we collect the definitions that are relevant to say what analytic 
functions between open subsets of locally convex spaces are. 

\begin{defn} \mlabel{def:0.2} (a) 
Let $E$ and $F$ be locally convex spaces, $U
\subeq E$ open and $f \: U \to F$ a map. Then the {\it derivative
  of $f$ at $x$ in the direction $h$} is defined as 
$$ \dd f(x)(h) := (\partial_h f)(x) := \lim_{t \to 0} \frac{1}{t}(f(x+th) -f(x)) $$
whenever it exists. The function $f$ is called {\it differentiable at
  $x$} if $\dd f(x)(h)$ exists for all $h \in E$. It is called {\it
  continuously differentiable}, if it is differentiable at all
points of $U$ and 
$$ \dd f \: U \times E \to F, \quad (x,h) \mapsto \dd f(x)(h) $$
is a continuous map. Note that this implies that the maps 
$\dd f(x)$ are linear (cf.\ \cite[Lemma~2.2.14]{GN15}). 
The map $f$ is called a {\it $\cC^k$-map}, $k \in \N \cup \{\infty\}$, 
if it is continuous, the iterated directional derivatives 
$$ \dd^{j}f(x)(h_1,\ldots, h_j)
:= (\partial_{h_j} \cdots \partial_{h_1}f)(x) $$
exist for all integers $1\leq j \leq k$, $x \in U$ and $h_1,\ldots, h_j \in E$, 
and all maps $\dd^j f \: U \times E^j \to F$ are continuous. 
As usual, $\cC^\infty$-maps are called {\it smooth}. 

 (b) If $E$ and $F$ are complex locally convex spaces, then $f$ is 
called {\it complex analytic} if it is continuous and, for each 
$x \in U$, there exists a $0$-neighborhood $V$ with $x + V \subeq U$ and 
continuous homogeneous polynomials $\beta_k \: E \to F$ of degree $k$ 
such that, for each $h \in V$, we have 
$f(x+h) = \sum_{k = 0}^\infty \beta_k(h),$
as a pointwise limit (\cite{BS71}). 
The map $f$ is called {\it holomorphic} if it is $\cC^1$ 
and, for each $x \in U$, the 
map $\dd f(x) \: E \to F$ is complex linear. 
If $F$ is sequentially complete, then $f$ is holomorphic if and only if 
it is complex analytic (\cite[Ths.~3.1, 6.4]{BS71}). 

(c) If $E$ and $F$ are real locally convex spaces, 
then we call a map $f \: U \to F$, $U \subeq E$ open, 
{\it real analytic} or a $\cC^\omega$-map, 
if for each point $x \in U$, there exists an open neighborhood 
$V \subeq E_\C$ and a holomorphic map $f_\C \: V \to F_\C$ with 
$f_\C\res_{U \cap V} = f\res_{U \cap V}$. 
The advantage of this definition, which differs from the one in 
\cite{BS71}, is that it also works nicely for non-complete spaces. 
Any analytic map is smooth, 
and the corresponding chain rule holds without any condition 
on the underlying spaces (see \cite{Gl02} for details).
\end{defn}

{\bf Acknowledgments:} 
D.~Belti\c t\u a was supported by the 
Grant of the Romanian National Authority for Scientific Research, CNCS-UEFISCDI, 
project number PN-II-ID-PCE-2011-3-0131, 
and the Emerging Fields Project ``Quantum Geometry'' of the University of Erlangen,  
and K.-H. Neeb was supported by DFG-grant NE 413/7-2, Schwerpunktprogramm ``Darstellungstheorie''.

\end{document}